\documentclass[a4paper, 11pt, final]{amsart}

\makeatletter

\numberwithin{equation}{section}
\numberwithin{figure}{section}
\usepackage[a4paper]{geometry}
\usepackage{amsmath,amsfonts,amssymb,amsthm,epsfig}

\usepackage{color,xcolor}
\usepackage[final,linkcolor = blue,citecolor = blue,colorlinks=true]{hyperref}
\usepackage[leqno]{amsmath}

\usepackage{amscd}
\usepackage{mathrsfs}
\usepackage{}
\usepackage[all]{xy}
\usepackage[OT1]{fontenc}
\usepackage[latin1]{inputenc}
\usepackage{amssymb,latexsym}
\usepackage{type1cm,ae}
\usepackage[notref,notcite]{showkeys}
\usepackage{graphicx}
\usepackage{xspace}
\usepackage{setspace}
\usepackage{enumerate}
\usepackage[english]{babel}
\usepackage{bbm}
\usepackage{esint}
\usepackage{soul}

\soulregister{\cite}7
\soulregister{\citep}7
\soulregister{\citet}7
\soulregister{\ref}7
\soulregister{\pageref}7
\sethlcolor{yellow}

\theoremstyle{definition}
\newtheorem{theorem}{Theorem}[section]
\newtheorem{proposition}[theorem]{Proposition}

\newtheorem{lemma}[theorem]{Lemma}

\newtheorem{corollary}[theorem]{Corollary}

\newtheorem{remark}[theorem]{Remark}

\numberwithin{equation}{section}

\newcommand*{\Wert}{\mathord{\mbox{|\kern-1.5pt|\kern-1.5pt|}}}

\DeclareMathOperator{\supp}{supp}

\def\rr{{\Bbb R}}
\def\rz{{{\rr}^n}}

\def\nn{{\Bbb N}}
\def\cc{{\Bbb C}}

\def\L{\mathcal{L}}

\def\H{\mathcal{H}}

\def\rdd{\rr_+^{n+1}}
\def\rdm{\rr^{n+1}}

\def\supp{{\rm{\ supp\ }}}

\def\ez{\epsilon}

\def\supp{{\rm supp}}

\def\l{\left}
\def\r{\right}

\def\XXint#1#2#3{{\setbox0=\hbox{$#1{#2#3}{\int}$}
		\vcenter{\hbox{$#2#3$}}\kern-.5\wd0}}

\makeatother

\usepackage{babel}

\begin{document}
\raggedbottom
\allowdisplaybreaks

\title[On regularity for parabolic systems with unbounded coefficients]{On regularity for nonhomogeneous parabolic systems with a skew-symmetric part in BMO}
	
	\author[Guoming Zhang]{Guoming Zhang}
	
	\address[Guoming Zhang]{College of Mathematics and System Science, Shandong University of Science and Technology, Qingdao, 266590, Shandong, People's Republic of China}
	\email{ zhangguoming256@163.com}

	

	\thanks{$^*$Corresponding author: Guoming Zhang}
	\thanks{The author is supported by the National Natural Science Foundation of Shandong Province (No. ~ZR2023QA124).}

	\date{}
	
	\begin{abstract} 
	In this paper we investigate the improved Caccioppoli inequality and the reverse H\"{o}lder inequality for gradients of weak solutions to nonhomogeneous parabolic systems whose coefficients can be split into a complex-valued and bounded part, which also satisfies the uniform G$\mathring{a}$rding inequality, and a real and anti-symmetric part in BMO. In particular, unbounded coefficients are allowed. 
	\end{abstract}

	\maketitle
	\tableofcontents
\section{Introduction}

In \cite{ABES}, Auscher, Bortz, Egert and Saari considered the local regularity of weak solutions to non-homogeneous parabolic systems of the form \begin{equation}\label{eq: 1.1}\partial_{t}u-\mbox{div}_{x}(A(t, x)\nabla_{x} u)=f+\mbox{div}_{x} F.\end{equation} Such a system is defined in an open parabolic cylinder $\Omega:=I_{0}\times Q_{0} \subset \rdm:=\rr\times \rz,$ where $I_{0}$ denotes an open interval in $\rr,$ and $Q_{0}$ is an open ball of $\rz.$ The coefficents are assumed to be elliptic in the  G$\mathring{a}$rding sense, bounded and measurable, notably, without any regularity on time as well as on all spatial variables. In their paper, they revealed a new regularity property (primarily referring to \cite[Theorem 8.1]{ABES}) of weak solutions to \eqref{eq: 1.1}, in contrast to the results in \cite{L, NS, GS}. This property can be roughly summarized as follows: 
\begin{theorem}\label{thm: 0.1}\; \emph{ If $u$ is a locally weak solution to \eqref{eq: 1.1} with $f$ only belonging to $L^{2_{*}}_{loc}(\Omega)$ ($2_{*}<2$) and $F\in L^{2}_{loc}(\Omega)$, then in time $u$ is bounded and H\"{o}lder continuous locally with values in spatial $L^{p}_{loc}(Q_{0})$ for some $p>2.$ } 
\end{theorem} As mentioned in \cite[Section 1.2]{ABES}, their results also generalize certain statements in \cite{C, GS}. 

Unlike the usual way, which relies on the local variational methods rooted in the Lions theory \cite{L} to construct the regularity properties of solutions to parabolic systems, the authors \cite{ABES} introduced a global variational approach. This methods utilizes the idea of splitting $\partial_{t}=D_{t}^{1/2}H_{t}D_{t}^{1/2}$ and explores the hidden coercivity within a parabolic energy space. This technique was first proposed by Kaplan \cite{K} and later studied for boundary value problems (BVPs) in \cite{HL, AEN, AK}. Equipped with such an approach, two distinct proofs for Theorem \ref{thm: 0.1} were presented: a real analysis proof in \cite[Section 6]{ABES}, based on a non-local Gehring's lemma, and an operator theoretic proof in \cite[Section 7]{ABES}, relying on the $\breve{S}$ne$\breve{I}$berg's theorem \cite{IS}.  

One of the objectives of this work is to show that the bounded assumption imposed on the coefficient matrix $A$ in Theorem \ref{thm: 0.1} can be removed. More precisely, we study the nonhomogeneous parabolic system with $$A(t, x)=S(t, x)+D(t, x),$$ where $S\in L^{\infty}(\rr^{n+1}, \cc^{mn})$ is complex-valued and satisfies the G$\mathring{a}$rding inequality (uniformly in $t$), i.e., for any $u\in W^{1, 2}(\rr^{n+1}, \cc^{m}),$ \begin{equation}\label{eq: 1.2}\begin{aligned}\mbox{Re} \;\int_{\rz}\;S(t,x)\nabla u(t, x)\cdot \overline{\nabla u(t, x)}&:=\mbox{Re} \;\int_{\rz}\;S^{\alpha \beta}_{kj}(t, x)\partial_{k}u^{\alpha}\overline{\partial_{j}u^{\beta}}\\
&\quad\quad\geq \lambda\int_{\rz}|\nabla u|^{2}-\kappa \int_{\rz}| u|^{2}\end{aligned}\end{equation} holds for some $\lambda, \kappa>0.$ The matrix $D$ is assumed to be real-valued, skew-symmetric and block diagonal, that is, $D=\l(\delta_{\alpha \beta}D^{\alpha \beta}\r)_{1\leq \alpha, \beta\leq m},$ with $D^{\alpha \beta}$ only belonging to the space $L^{\infty}(\rr; BMO(\rz)).$ Hence, $A$ is unbounded. Under these hypotheses on $A,$ we state our main theorem as follows. 

\begin{theorem}\label{thm: 4.1}\;  Let $f\in L_{loc}^{2_{*}}(\Omega, \cc^{m})$ and $F\in L^{2}_{loc}(\Omega, \cc^{mn}).$  Suppose that $v$ is a weak solution of $\partial_{t}v-\mbox{div}((S(t, x)+D(t, x))\nabla v)=f+\mbox{div} F\; \mbox{in}\;\Omega.$ Let $\gamma>1$ and $I\times Q$ be a open paraboli cylinder with $l(I)\approx r(Q)^{2}$ such that $\overline{\gamma^{2}I\times \gamma Q}\subset \Omega.$ If $p>2$ is sufficiently close to 2 and $\alpha:=\frac{1}{2}-\frac{1}{p},$ then 
\begin{equation*}\begin{aligned}
\l(\fint_{I\times Q}|\nabla v|^{p}\r)^{1/p}&+\sup_{t\in I}\l(\fint_{I\times Q}| v(t, \cdot)|^{p}\r)^{1/p}+\sup_{t, s\in I}\l(\fint_{Q}\frac{|v(t, \cdot)-v(s, \cdot)|^{p}}{|t-s|^{\alpha p}}\r)^{1/p}\\
&\lesssim \frac{1}{r(Q)}\l(\fint_{\gamma^{2}I\times \gamma Q}| v|^{2}\r)^{1/2}+\l(\fint_{\gamma^{2}I\times \gamma Q}|F|^{p}\r)^{1/p}+r(Q)\l(\fint_{\gamma^{2}I\times \gamma Q}|f|^{p_{*}}\r)^{1/p_{*}}.
\end{aligned} \end{equation*}
The implicit constant depends only on dimensions, $\|D\|_{L^{\infty}(\rr; BMO(\rz))}$, $\gamma,$ ellipticity and the constants controlling the ratio $l(I)/r(Q)^{2}.$
\end{theorem}

\begin{theorem}\label{thm: 4.2}\;  Under the same hypotheses as in Theorem \ref{thm: 4.1},  we have  
\begin{equation}\label{eq: 4.1}
\l(\fint_{I\times Q}|\nabla v|^{p}\r)^{1/p}\lesssim \fint_{\gamma^{2}I\times \gamma Q}| \nabla v|+\l(\fint_{\gamma^{2}I\times \gamma Q}|F|^{p}\r)^{1/p}+r(Q)\l(\fint_{\gamma^{2}I\times \gamma Q}|f|^{p_{*}}\r)^{1/p_{*}},\end{equation}
where the number $p$ is greater than 2 but sufficiently close to 2.  The implicit constant depends only on $\|D\|_{L^{\infty}(\rr; BMO(\rz))}$, dimensions, $\gamma,$ ellipticity and the constants controlling the ratio $l(I)/r(Q)^{2}.$
\end{theorem} 

Theorem \ref{thm: 4.1} and Theorem \ref{thm: 4.2} serves as natural extensions of \cite[Theorem 8.1]{ABES} and \cite[Theorem 8.2]{ABES}, respectively.

Our motivation stems from the investigation of the behavior of solutions of elliptic and parabolic equations with a divergence-free drift \cite{MV, LZ, O, ZQ, SS}. Most notably, the authors \cite{SS} studied equations of the form $\partial_{t}-\Delta u+\textbf{b}\cdot \nabla u=0$ resp. $-\Delta u+\textbf{b}\cdot \nabla u=0$ with $\mbox{div}\textbf{b}=0.$ They discovered that the usual regularity assumptions on $\textbf{b}$ can be relaxed by the divergence-free condition, and the classical Liouville theorem and Harnack inequality remain valid. It also turns out that the interior theory of De Giorgi, Nash and Moser is preserved for these parabolic operators if $\textbf{b}\in L^{\infty}(BMO^{-1})$ resp. elliptic operators if $\textbf{b}\in BMO^{-1}.$ For gerneral elliptic and parabolic opertaors with real coefficients in divergence form, this condition is equivalent to the coefficients matrix $A$ being expressible as the sum of an  $L^{\infty}$ elliptic symmetric part and an unbounded anti-symmetric part which belongs to $L^{\infty}(BMO)$ in parabolic case resp. $\textbf{b}\in BMO$ in elliptic case. The space BMO has a key role in two ways which will be explicitly displayed in following sections as well. First, its right scaling properties enable one to implement the iterative arguments of De Giorgi-Nash-Moser. Second, the BMO condition on the anti-symmetric part of the matrix, combined with the compensated compactness as discussed in \cite{CLMS, SS, LP}, allows for the proper definition of weak solutions.  

Operators with a real and skew-symmetric part in $BMO,$ of the form $\L:=-\mbox{div}(A\nabla )$ or $\H:=\partial_{t}-\mbox{div}(A\nabla ),$ have gained much attention since \cite{SS}. In \cite{XQ}, Qian and Xi showed the existence of the fundamental solution for such parabolic operator $\H$ with real coefficients and derived the pointwise Gaussian bounds. The BVPs generated by such an elliptic operator $\L,$ where the symmetric part in $A$ is real, bounded and uniform elliptic, was launched in \cite{LP}. In this work, the authors constructed the Green function and elliptic measure associated to $\L,$ and proved the boundary H\"{o}lder continuity of the solution, the solvability of the continuous Dirichlet problem on non-tangentially accessible (NTA) domains and the well-posedness of the $L^{p}$ Dirichlet problem with respect to the elliptic measure. Later, Hofmann, Li, Mayboroda and Pipher \cite{HLMPJ} established the solvability of the Dirichlet problem generated by the same operator $\L$ with $t$-independent coefficients in upper half space $\rdd$ for some $p\in (1, \infty),$ equivalently, the elliptic measure associated to $\L$ belongs to the $A_{\infty}$ class with respect to the Lebesgue measure. In another direction, for $\L$ having a complex-valued, uniform elliptic and symmetric part, Escauriaza and Hofmann \cite{EH} addressed the $L^{2}$ Kato square root problem without using Gaussian estimate. Subsequently, adding the real-valued assumption on $A,$ Hofmann, Li, Mayboroda and Pipher \cite{HLMP} generalized these results to the $L^{p}$ setting. Regarding the Kato square root problem associated to $\H$ with $A=S+D,$ where $S\in L^{\infty}$ is complex-valued and uniform elliptic, and $D\in L^{\infty}(BMO)$ is real and anti-symmetric, Ataei and Nystr\"{o}m \cite{AK} provided a positive answer in $L^{2}$ case.


The aforementioned results are centrally concerned with generalizations to BMO skew-symmetric part for elliptic and parabolic equations. There is not much work on such an extension to elliptic and parabolic systems.  It is important to emphasize that Dong and Kim \cite{DK} extended the results on fundamental solutions presented in \cite{XQ} to second-order parabolic systems with real coefficients, under the assumption that weak solutions of these systems satisfy a specific local boundedness estimate. This advancement has been extended to scenarios involving complex-valued coefficients in the work of Auscher and Egert \cite{AE}. To our knowlege, this represents the most recent progress in the study of parabolic sysytems featuring a BMO skew-symmetric part. Therefore, our results have provided an update in this regard. 

In the next Section 2 we introduce relevant notation. Section 3 is devoted to the global variational setup associated to parabolic operators with a BMO skew-symmetric part. In Section 4 the improved Caccioppoli inequality is established. Building on this result, we provide the proofs of Theorem \ref {thm: 4.1} and Theorem \ref {thm: 4.2} in Section 5. An enhanced version of the Gehring lemma is included in Section 6 as an appendix.

\section{Preliminaries and basic notations}

For simplicity, we abbreviate $\nabla:=\nabla_{x}$ and $\mbox{div}:=\mbox{div}_{x}$ with respect to the gradient and divergence in the spatial variables $x$, respectively. For any $1<p<\infty,$ $p^{*}$ and $p_{*},$ defined by $$\frac{1}{p^{*}}:=\frac{1}{p}-\frac{1}{n+2},\quad \frac{1}{p_{*}}:=\frac{1}{p}+\frac{1}{n+2},$$ are separately referred to as the Sobolev upper and lower conjugate. As usual, $p'$ denotes the H\"{o}lder conjugate of $p$, then $2^{*}=2_{*}'.$

\noindent{\textbf{2.1. \;Weak solutions.}} \quad Given $f\in L^{1}_{loc}(\Omega; \cc^{m})$ and $F\in L^{1}_{loc}(\Omega; \cc^{mn}),$ we say that $u$ is a weak solution to $$\partial_{t}u-\mbox{div}((S(t, x)+D(t, x))\nabla u)=f+\mbox{div} F\quad \mbox{in}\; \Omega$$ if $u\in L^{2}_{loc}(I; W^{1, 2}_{loc}(Q; \cc^{m}))$ satisfies  \begin{equation}\label{eq: 1.3}\int_{\Omega}(S+D)\nabla u\cdot \overline{\nabla \varphi}-u\overline{\partial_{t}\varphi}=\int_{\Omega}f\cdot \overline{\varphi}-F\cdot \overline{\nabla \varphi}\quad\mbox{for any}\; \varphi\in C_{0}^{\infty}(\Omega).\end{equation}

\noindent{\textbf{2.2. \;Function spaces.}} \quad The operators $D_{t}^{1/2}$ and $H_{t}$ represent the half-order time derivative and Hilbert transform in time through the Fourier symbols $|\tau|^{1/2}$ and $isgn(\tau),$ respectively. We define the space $H^{1/2}(\rr; L^{2}(\rr^{n}))$ as the set of functions $f\in L^{2}(\rdm)$ such that $D_{t}^{1/2}f\in L^{2}(\rdm).$ Obviously, the space of smooth and compactly supported functions $C_{0}^{\infty}(\rdm)$ serves as a dense subspace of $H^{1/2}(\rr; L^{2}(\rr^{n})).$ The function space $V:=L^{2}(\rr; W^{1, 2}(\rz)),$ endowed with $$\|u\|_{V}:=\l( \|u\|^{2}_{L^{2}}+\|\nabla u\|^{2}_{L^{2}}\r)^{1/2},$$ is a Hilbert space, which gives rise to a natural space for global weak solutions to parabolic equations. Furthermore, we also need the Hilbert space $E:=V\cap H^{1/2}(\rr; L^{2}(\rr^{n})):=L^{2}(\rr; W^{1,2}(\rr^{n}))\cap H^{1/2}(\rr; L^{2}(\rr^{n}))$ with norm defined by $$\|u\|_{E}:=\l( \|u\|_{L^{2}}^{2}+ \|\nabla u\|^{2}_{L^{2}}+\|D_{t}^{1/2} u\|^{2}_{L^{2}}\r)^{1/2}.$$

\noindent{\textbf{2.3. \;Floating constants.}}\quad For ease of presentation,  we adopt the notations $\lesssim$ and $\approx.$ Specifically, for two positive constants $a, b$, the expression $a\lesssim b$ means that there exists a nonessential constant $C$ such that $a\leq C b$, while $a\approx b$ signifies that both $a\leq C b$ and $b\leq C a$ hold. The value of the constant $C$ might change at each occurrence which depends upon $n, \lambda, \kappa$ and other nonessential constants.

\section{The global variational setup} 	
	
We start this section by referencing two key propositions. The first one is derived from \cite{LP, SS} and is referred to as the compensated compactness estimate.. The second one proved in \cite{HLMP} is a variant of it. These are indispensable tools for addressing the contributions of the anti-symmetric BMO part in $A$. 
 
\begin{proposition}\label{pro: 2.1}(\cite{LP, SS})\; Let $1<p<\infty,$ $u\in \dot{W}^{1, p}(\rz)$ and $v\in \dot{W}^{1, p'}(\rz).$ Then for any $1\leq k, j\leq n,$ $\partial_{k}u\partial_{j}v-\partial_{j}u\partial_{k}v$ belongs to Hardy space $ \mathcal{H}^{1}(\rz)$ and $$\|\partial_{k}u\partial_{j}v-\partial_{j}u\partial_{k}v\|_{\mathcal{H}^{1}(\rz)}\lesssim \|\nabla u\|_{L^{p}(\rz)}\|\nabla v\|_{L^{p'}(\rz)}.$$
\end{proposition}

\begin{proposition}\label{pro: 2.2}(\cite[Proposition 2.6]{HLMP})\;  Let $u, v\in W^{1, 2}(\rz)$ and $\varphi$ be a Lipschitz function on $\rz.$  Then for any $1\leq k, j\leq n,$ $\partial_{k}(uv)\partial_{j}\varphi-\partial_{j}(uv)\partial_{k}\varphi\in \mathcal{H}^{1}(\rz)$ with $$\|\partial_{k}(uv)\partial_{j}\varphi-\partial_{j}(uv)\partial_{k}\varphi\|_{\mathcal{H}^{1}(\rz)}\lesssim\|\nabla \varphi\|_{L^{\infty}}(\|u\|_{L^{2}(\rz)}\|\nabla v\|_{L^{2}(\rz)}+\|v\|_{L^{2}(\rz)}\|\nabla u\|_{L^{2}(\rz)}).$$
\end{proposition}

In the sequel, we sometimes use the notation for the scalar case $m =1$ (even
if $m > 1$) when we are only interested in norm estimates.
Lemma \ref{lemma: 2.3} below states that the operator $$\H_{\kappa}:=\partial_{t}-\mbox{div}((S(t, x)+D(t, x))\nabla )+\kappa+1$$ is invertible from $E$ to its dual space $E^{*},$ which is a starting point of our proof. In particular, the case $\kappa=0$ is addressed in \cite[Lemma 4.1]{AK}. 
	
\begin{lemma}\label{lemma: 2.3}\;  Assume \begin{equation}\label{eq: 2.1}<\H_{\kappa} u, v>:=\int_{\rr^{n+1}}A\nabla u\cdot \overline{\nabla v}+H_{t}D_{t}^{1/2}u\cdot \overline{D_{t}^{1/2}v}+(\kappa+1)u\cdot\overline{v},\quad u, v\in E.\end{equation} Then $\H$ is bounded from $E$ to $E^{*}$ and its inverse is also bounded with both norms depending only on ellipticity, $\|D\|_{L^{\infty}(\rr; BMO(\rz))}$
and dimensions.\end{lemma} 
{\it Proof.}\quad By the assumption on $D$ and Proposition \ref{pro: 2.1}, 
we see \footnote{For the validity of the equation involving the integration of $D$ one can see \cite[Section 2.16]{AE}.}   \begin{equation}\label{eq: 2.2}
\begin{aligned}
\int_{\rr^{n+1}}D\nabla u\cdot \overline{\nabla v}
&=\frac{1}{2} \sum_{\alpha=\beta}\int_{\rr^{n+1}}D_{kj}^{\alpha \beta}(\partial_{k}u^{\alpha}\overline{\partial_{j}v^{\beta}}-\partial_{j}u^{\alpha}\overline{\partial_{k}v^{\beta}}) \\
&\lesssim \int_{\rr}\|D(t, \cdot)\|_{BMO(\rz)}\|\nabla u\|_{L^{2}(\rz)}\|\nabla v\|_{L^{2}(\rz)} \\
&\lesssim \|\|D(t, \cdot)\|_{BMO(\rz)}\|_{L^{\infty}(\rr)}\|\nabla u\|_{L^{2}(\rr^{n+1})}\|\nabla v\|_{L^{2}(\rr^{n+1})}.\end{aligned} 
\end{equation} So $\H$ is bounded from $E$ to $E^{*}.$ Following \cite[Lemma 3.2]{ABES}, we define the sesquilinear form on $E\times E$ by
\begin{equation*}a_{\delta}(u, v):=\int_{\rr^{n+1}}A\nabla u\cdot \overline{\nabla (1+\delta H_{t})v}+H_{t}D_{t}^{1/2}u\cdot \overline{D_{t}^{1/2}(1+\delta H_{t})v}+(\kappa+1)u\cdot\overline{(1+\delta H_{t})v}.\end{equation*} As in \eqref{eq: 2.2}, using also the boundedness of $H_{t}$ on $L^{2},$ it follows that $a_{\delta}(\cdot, \cdot)$ is bounded. To prove the invertibility, we first claim that $a_{\delta}$ is accretive for sufficiently small enough $\delta$. Using $S\in L^{\infty}$, \eqref{eq: 2.1}, $\int_{\rr^{n+1}} \mbox{Re}\; H_{t}D_{t}^{1/2}u\cdot \overline{D_{t}^{1/2}u}=0$ and the fact that $H_{t}$ commutes with $D_{t}^{1/2}$ and $\nabla,$ we obtain 
\begin{equation*}\begin{aligned} 
\mbox{Re}\;a_{\delta}(u, u)&\geq (\lambda-\|S\|_{L^{\infty}})\|\nabla u\|_{L^{2}(\rr^{n+1})}^{2}+\delta \|D_{t}^{1/2} u\|_{L^{2}(\rr^{n+1})}^{2}+\|u\|_{L^{2}(\rr^{n+1})}^{2}\\
&\quad -\mbox{Re}\;\int_{\rr^{n+1}}D\nabla u\cdot \overline{\nabla (1+\delta H_{t})u}.\end{aligned} \end{equation*}
We note that by the anti-symmetry of $D,$ \begin{equation}\label{eq: 2.3}\mbox{Re}\;D\xi\cdot \overline{\xi}=\frac{1}{2}\sum_{\alpha=\beta}D_{kj}^{\alpha \beta}(\xi_{k}\overline{\xi_{j}}-\xi_{j}\overline{\xi_{k}})=0\;\;\mbox{for any}\;\; \xi\in \cc^{mn},\end{equation} moreover, by \eqref{eq: 2.2}, $$\bigg|\int_{\rr^{n+1}}D\nabla u\cdot \overline{\nabla \delta H_{t}u}\bigg|\lesssim\delta \|\|D(t, \cdot)\|_{BMO(\rz)}\|_{L^{\infty}(\rr)}\|\nabla u\|^{2}_{L^{2}(\rr^{n+1})}.$$ The claim is proved. Then the Lax-Milgram theorem yields that $$<\H_{\kappa} u, (1+\delta H_{t})v>=a_{\delta}(u, (1+\delta H_{t})v)\;\;\mbox{for any}\;\;u, v\in E.$$ This immediately implies the invertibility of $\H_{\kappa}$ from $E$ to $E^{*},$ since $(1+\delta H_{t})$ is isometric on $E$ when $\delta$ is small enough.

\hfill$\Box$ 	
	
To proceed, we require two lemmas in \cite{ABES}. Set $V_{p}:=V\cap L^{p}(\rdm)$ $(1<p<\infty),$ equipped with the norm $\|u\|_{V_{p}}:=\max\{\|u\|_{V}, \|u\|_{L^{p}}\}.$ Then, clearly, $V_{2}=V.$ We also let $V_{p}^{*}$ denote its dual space. In light of \cite[Theorem 2.7.1]{BL}, the $V_{p}^{*}-V_{p}$ duality can be realized as a Lebesgue integral (see the formular over \cite[lemma 3.3]{ABES}), which is used there to deduce the following lemma:

\begin{lemma}\label{lemma: 2.41}(\cite[lemma 3.3]{ABES})\;  Consider a function in $v\in V_{p}$ such that $\partial_{t}v \in V^{*}_{p}.$ Then $u\in C(\rr; L^{2}(\rz))$ and $t\to \|v(t, \cdot)\|_{2}^{2}$ is absolutely continous on $\rr,$ vanishes at $\pm \infty,$ and satisfies $$\sup_{t\in \rr}\|v(t, \cdot)\|_{2}^{2} \lesssim \|v\|_{V_{p}} \|\partial_{t}v\|_{V^{*}_{p}}.$$ Moreover, it holds $\mbox{Re}<\partial_{t} v, v>=0$ for the $V_{p}^{*}-V_{p}$ duality.\end{lemma}

Lemma \ref{lemma: 2.42} below essentially establishes the $L^{p}(\rdm)\to L^{p^{*}}(\rdm)$ boundedness of the parabolic Reisz potential, as discussed in \cite{GR}.
\begin{lemma}\label{lemma: 2.42}(\cite[lemma 3.4]{ABES})\;  Let $1<p<n+2.$ Then for any $\phi\in C_{0}^{\infty}(\rdm),$ $$\|\phi\|_{L^{p^{*}}}\lesssim \|\nabla \phi\|_{L^{p}} +\|D^{1/2}_{t}\phi\|_{L^{p}}.$$
\end{lemma}

The next proposition extends \cite[Proposition 3.1]{ABES} and ensures the higher integrability of the solution to \eqref{eq: 1.1}, palying a cnetral role in our proof.

\begin{proposition}\label{pro: 2.4} Let $f\in L^{2_{*}}(\rr^{n+1})$ and $F\in L^{2}(\rr^{n+1}).$ If $v\in V$ is a weak solution of $\partial_{t}v-\mbox{div}((S(t, x)+D(t, x))\nabla v)=f+\mbox{div} F$ in $\rdm,$ then 
\begin{equation*}
\begin{aligned} 
&(i)\quad v\in  H^{1/2}(\rr, L^{2}(\rr^{n}));\\
&(ii)\quad v\in L^{2^{*}}(\rr^{n+1});\\
&(iii)\quad v\in C_{0}(\rr, L^{2}(\rz)) \;\mbox{and}\;\;\|v(t, \cdot)\|^{2}_{L^{2}(\rz)}\; \mbox{is absolutely continuous in} \;\; t\in \rr,
\end{aligned} 
\end{equation*}
moreover, the inequality \begin{equation}\label{eq: 2.4}\sup_{t\in \rr}\|v(t, \cdot)\|_{L^{2}(\rz)}+\|v\|_{L^{2^{*}}(\rr^{n+1})}+\|D_{t}^{1/2}v\|_{L^{2}(\rdm)}\lesssim\|v\|_{V}+\|f\|_{L^{2_{*}}(\rr^{n+1})}+\|F\|_{L^{2}(\rdm)}\end{equation} holds. The implicit constant depends only on dimensions, $\|D\|_{L^{\infty}(\rr; BMO(\rz))}$ and ellipticity.
\end{proposition}

{\it Proof.}\quad By utilizing self-evident embeddings and repeating the arguments from the beginning of the proof of \cite[Proposition 3.1]{ABES} we can show $f+(\kappa+1)v+\mbox{div} F \in E^{*}.$ Thus it follows from Lemma \ref{lemma: 2.3} that there exists a unique $\tilde{v}\in E$ such that $\H_{\kappa} \tilde{v}=f+(\kappa+1)v+\mbox{div} F.$ This indicates that, for any $\psi\in E,$ $$\int_{\rr^{n+1}}(S+D)\nabla \tilde{v}\cdot \overline{\psi}+H_{t}D_{t}^{1/2}\tilde{v}\cdot \overline{D_{t}^{1/2}\psi}+(\kappa+1)\tilde{v}\cdot\overline{D_{t}^{1/2}\psi}=\int_{\rr^{n+1}}(f+(\kappa+1)v)\overline{\psi}-F\cdot \overline{\nabla \psi}.$$ By restricting to $\psi\in C_{0}^{\infty}(\rdm)$ we can write $$\int_{\rdm}H_{t}D_{t}^{1/2}\tilde{v}\cdot \overline{D_{t}^{1/2}\psi}=-\int_{\rdm}\tilde{v}\overline{\partial_{t}\psi}.$$ As a result, by \eqref{eq: 1.3}, $u:=v-\tilde{v} \in V$ is a weak solution of \begin{equation}\label{eq: 2.5}\partial_{t}u-\mbox{div}((S(t, x)+D(t, x))\nabla u)+(\kappa+1)u=0 \;\;\mbox{in}\; \rdm.\end{equation} In the current setting, we assert that: \begin{equation}\label{eq: 2.6}\fint_{I\times Q}|\nabla u|^{2}\leq C(\ez) \frac{1}{r(Q)^{2}}\fint_{4I\times 2Q}| u|^{2}+\ez \fint_{4I\times 2Q}| \nabla u|^{2}\end{equation} for any parabolic cylinder $I\times Q$ with $l(I)= r(Q)^{2},$ where $l(I), r(Q)$ denote the length of $I$ and the radius of $Q.$ \eqref{eq: 2.6} is proved via a standard truncation method, as outlined in Theorem \ref{thm: 3.1}, a technique that will be repeatedly employed in subsequent sections.

By the scaling invariance property of BMO functions, we can assume, without loss of generality, that $r(Q)=1.$ Let $\chi, \eta \in C_{0}^{\infty}(\rdm)$ be two real-valued functions satisfying the following conditions: 

\emph{ $\chi\equiv 1$ on $I\times Q$ with its support contained in $\frac{9}{4}I\times \frac{3}{2}Q,$ and $\eta \equiv 1$ on $\frac{9}{4}I\times \frac{3}{2}Q$ with its support contained in $4I\times 2Q.$}\\
 Also, we let $v:=u\chi.$ According to Theorem \ref{thm: 3.1}, $v$ is a weak solution to $$\partial_{t}v-\mbox{div}\l(\l(S(t, x)+\l(D-\fint_{ 2Q}D\r)\r)\nabla v\r)=\tilde{f}+\mbox{div} \tilde{F}\quad\mbox{in}\;\rdm$$ with $$\tilde{f}:=(\partial_{t}\chi)u-\l(S+D-\fint_{ 2Q}D\r)\nabla u\cdot \nabla \chi, \quad \tilde{F}:=-\l(S(t, x)+D-\fint_{4I\times 2Q}D\r)(u\nabla \chi).$$ 

Recall that $u\in V.$ Then $\partial_{t}v\in V^{*}$ (the dual space of $V$). In addition, by setting $$T_{S, D}:=\mbox{div}\l(\l(S(t, x)+\l(D-\fint_{ 2Q}D\r)\r)\nabla v\r),$$ we can show that $T_{S, D}, \tilde{f},$ and $\mbox{div} \tilde{F}$ also belong to $ V^{*}.$ For $\tilde{f},$ it suffices to note, by H\"{o}lder inequality and the $L^{2}(\rz) \to L^{\tilde{2}}(\rz)$ boundedness of the Reisz potential $I_{1}$ of order 1 on $\rz,$ \footnote{The requirement $n\geq 3$ is essential for the strong $(2, \tilde{2})$ boundedness of Reisz potential $I_{1}$ to hold. When $n=2,$ the limiting Sobolev inequality from $L^{\infty}$ to the Lorentz space $L^{(2, 1)}$ is expressed as $$\|\psi\|_{L^{\infty}}\lesssim \|\nabla \psi\|_{L^{(2, 1)}}, \quad \mbox{for any}\; \psi\in C_{0}^{\infty}.$$ Unfortunately, this is not strong enough to bound $T_{1}$ and $T_{2}$ as needed.} that 
\begin{align*}
&T_{1}:=\bigg|\int_{\rdm} \l(D-\fint_{ 2Q}D\r)\nabla u\cdot \chi\overline{\phi}\bigg|\\[4pt]
&\quad\quad\lesssim \int_{4I}\l(\int_{2Q}\bigg|D-\fint_{ 2Q}D\bigg|^{n}\r)^{1/n}\l(\int_{2Q}|\nabla u \chi|^{2}\r)^{1/2}\l(\int_{2Q}|\phi(t, x)|^{\tilde{2}}dx\r)^{1/\tilde{2}}\quad \l(\frac{1}{\tilde{2}}:=\frac{1}{2}-\frac{1}{n}\r)
\\[4pt]
&\quad\quad\lesssim \|\|D(t, \cdot)\|_{BMO(\rz)}\|_{L^{\infty}(\rr)}\int_{4I}\l(\int_{2Q}|\nabla u \chi|^{2}\r)^{1/2}\l(\int_{\rz}|\nabla \phi(t, x)|^{2}dx\r)^{1/2}\\[4pt]
&\quad\quad\lesssim \|\|D(t, \cdot)\|_{BMO(\rz)}\|_{L^{\infty}(\rr)}\|\nabla u\|_{L^{2}(\rdm)}\|\nabla \phi\|_{L^{2}(\rdm)}\end{align*} holds for any $\phi\in V.$ In a same manner, we derive the estimate: $$T_{2}:=\bigg|\int_{\rdm} \l(D-\fint_{ 2Q}D\r)(u\nabla \chi)\overline{\phi}\bigg|\lesssim \|\|D(t, \cdot)\|_{BMO(\rz)}\|_{L^{\infty}(\rr)}\|u\|_{L^{2}(\rdm)}\|\nabla \phi\|_{L^{2}(\rdm)},$$ which implies $\tilde{F}\in V^{*}.$ On the other hand, by Proposition \ref{pro: 2.1}, we have that $$\bigg|\int_{\rdm} \l(D-\fint_{ 2Q}D\r)\nabla v\cdot\overline{\nabla \phi}\bigg|\lesssim \|\|D(t, \cdot)\|_{BMO(\rz)}\|_{L^{\infty}(\rr)} \|\nabla v\|_{L^{2}(\rdm)}\|\nabla \phi\|_{L^{2}(\rdm)}.$$ Thus, $T_{S, D}$ belongs to $V^{*}.$ The above estimates in conjunction with Lemma \ref{lemma: 2.41} yield 
\begin{equation}\label{eq: 2.7}\begin{aligned}
 0&=-\mbox{Re}\;\int_{\rdm}S\nabla v\cdot \overline{\nabla v}-S(u\nabla \chi)\cdot \overline{\nabla v}+\mbox{Re}\;\int_{\rdm}\l(D-\fint_{2Q}D\r)(u\nabla \chi)\cdot \overline{\nabla v}\\
 &\quad +\mbox{Re}\;\int_{\rdm}\partial_{t}\chi u\overline{ v}-\mbox{Re}\;\int_{\rdm}S \nabla u\cdot \nabla\chi\overline{ v}-\mbox{Re}\;\int_{\rdm}\l(D-\fint_{2Q}D\r)\nabla u\cdot \nabla \chi\overline{ v}\\
&\quad\quad+\mbox{Re}\;\int_{\rdm} (\kappa+1)v\cdot \overline{ v}\\
&:=-U_{1}+U_{2}+U_{3}-U_{4}-U_{5}-U_{6}.
\end{aligned} \end{equation}
Due to the assumption $S\in L^{\infty}(\rdm),$ we may follow, mutatis mutandis, the proof of \cite[Proposition 4.3]{ABES} or Proposition \ref{pro: 3.2} below, to manage $U_{1}, U_{3}, U_{4},$ and $U_{6}.$ At this point, it is enough to exploit H\"{o}lder inequality and Young's inequality to infer that \begin{equation}\label{eq: 2.8}|U_{1}+U_{3}+U_{4}|\lesssim \|u\|^{2}_{L^{2}(4 I\times 2Q)}+\ez \|\nabla u\|^{2}_{L^{2}(4 I\times 2Q)}+\delta \|\nabla v\|^{2}_{L^{2}(\rdm)},\end{equation} where $\delta, \ez>0$ are small enough. Hence, invoking the G$\mathring{a}$rding inequality \eqref{eq: 1.2} and \eqref{eq: 2.7}-\eqref{eq: 2.8}, and taking the limit as $\delta\to 0,$  we get \begin{equation}\label{eq: 2.9}\|\nabla v\|^{2}_{L^{2}(\rdm)}\leq C(\ez)\|u\|^{2}_{L^{2}(4 I\times 2Q)}+\ez \|\nabla u\|^{2}_{L^{2}(4 I\times 2Q)}+|U_{2}|+|U_{5}|.\end{equation}

We now discuss the approach to handle $U_{2}$ and $U_{5}.$ Clearly, the two terms are of a similar nature because of $v=u\chi.$ For this reason, we only illustrate the procedure for dealing with $U_{5}.$ To the end, we first set $$\tilde{D}:=D-\fint_{2Q}D=\l(D_{kj}^{\alpha \beta}-\fint_{2Q}D_{kj}^{\alpha \beta}\r).$$ Then $\tilde{D}$ is skew-symmetric and $\|\|\tilde{D}(t, \cdot)\|_{BMO(\rz)}\|_{L^{\infty}(\rr)}=\|\|D(t, \cdot)\|_{BMO(\rz)}\|_{L^{\infty}(\rr)}.$ Using \begin{equation*}\begin{aligned}
U_{5}&=\frac{1}{2}\mbox{Re}\;\sum_{\alpha=\beta}\int_{\rr^{n+1}}\tilde{D}_{kj}^{\alpha \beta}\partial_{k}u^{\alpha}\partial_{j}\chi^{2}\overline{u^{\beta}}\\
&=\frac{1}{2}\sum_{\alpha=\beta}\mbox{Re}\;\int_{\rr^{n+1}}\tilde{D}_{kj}^{\alpha \beta}\partial_{k}(u^{\alpha}\overline{u^{\beta}})\partial_{j}\chi^{2}\\
&\quad-\frac{1}{2}\sum_{\alpha=\beta}\mbox{Re}\;\int_{\rr^{n+1}}\tilde{D}_{kj}^{\alpha \beta}\partial_{k}\overline{u^{\beta}}u^{\alpha}\partial_{j}\chi^{2}\end{aligned}\end{equation*} and $$\sum_{\alpha=\beta}\mbox{Re}\;\int_{\rr^{n+1}}\tilde{D}_{kj}^{\alpha \beta}\partial_{k}\overline{u^{\beta}}u^{\alpha}\partial_{j}\chi^{2}=\mbox{Re}\;\sum_{\alpha=\beta}\int_{\rr^{n+1}}\tilde{D}_{kj}^{\alpha \beta}\partial_{k}u^{\alpha}\partial_{j}\chi^{2}\overline{u^{\beta}},$$ together with the anti-symmetry of $D^{\alpha \beta}$ and the fact that $\eta\equiv 1$ on $\supp \chi,$ we conclude that \begin{equation*}\begin{aligned}U_{5}&=\frac{1}{4}\sum_{\alpha=\beta}\mbox{Re}\;\int_{\rr^{n+1}}\tilde{D}_{kj}^{\alpha \beta}\partial_{k}(u^{\alpha}\overline{u^{\beta}})\partial_{j}\chi^{2}\\
&=\frac{1}{8}\sum_{\alpha=\beta}\mbox{Re}\;\int_{\rr^{n+1}}\tilde{D}_{kj}^{\alpha \beta}(\partial_{k}(u^{\alpha}\overline{u^{\beta}})\partial_{j}\chi^{2}-\partial_{k}\chi^{2}\partial_{j}(u^{\alpha}\overline{u^{\beta}}))\\
&=\frac{1}{8}\sum_{\alpha=\beta}\mbox{Re}\;\int_{\rr^{n+1}}\tilde{D}_{kj}^{\alpha \beta}(\partial_{k}(u^{\alpha}\overline{u^{\beta}}\eta^{2})\partial_{j}\chi^{2}-\partial_{k}\chi^{2}\partial_{j}(\eta^{2}u^{\alpha}\overline{u^{\beta}})). \end{aligned} \end{equation*} With the latter equality in hand, an application of Proposition \ref{pro: 2.2} gives \begin{equation*}\begin{aligned}|U_{5}|&\lesssim\int_{4 I}\|D_{kj}^{\alpha \beta}(t, \cdot)\|_{BMO(\rz)}\|(\partial_{k}(u^{\alpha}\overline{u^{\beta}}\eta^{2})\partial_{j}\chi^{2}-\partial_{k}\chi^{2}\partial_{j}(\eta^{2}u^{\alpha}\overline{u^{\beta}}))\|_{\mathcal{H}^{1}(\rz)}\\
&\lesssim\int_{ 4 I}\|D(t, \cdot)\|_{BMO(\rz)}|\|\nabla\chi^{2} \|_{L^{\infty}}\|u\eta\|_{L^{2}(\rz)}\|\nabla(u\eta)\|_{L^{2}(\rz)}\\
&\lesssim\int_{4 I}\|D(t, \cdot)\|_{BMO(\rz)}|\|\nabla\chi^{2}\|_{L^{\infty}} \l(\|u\eta\|_{L^{2}(\rz)}\|u\nabla\eta\|_{L^{2}(\rz)}+\|u\eta\|_{L^{2}(\rz)}\|\nabla u\eta\|_{L^{2}(\rz)}\r)\\
&\lesssim \|\|D(t, \cdot)\|_{BMO(\rz)}\|_{L^{\infty}(\rr)}\l(\int_{ 4I }\|u\eta\|^{2}_{L^{2}(\rz)}\r)^{1/2}\l(\int_{4 I}\|u\nabla\eta\|^{2}_{L^{2}(\rz)}\r)^{1/2}
\end{aligned} \end{equation*}
\begin{equation*}\begin{aligned}
&\quad\quad+ \|\|D(t, \cdot)\|_{BMO(\rz)}\|_{L^{\infty}(\rr)}\l(\int_{ 4 I}\|u \eta\|^{2}_{L^{2}(\rz)}\r)^{1/2}\l(\int_{ 4 I}\|\nabla u\eta \|^{2}_{L^{2}(\rz)}\r)^{1/2}\\
&:=U_{51}+U_{52}.
\end{aligned} \end{equation*}
Obviously, $$U_{51} \lesssim \|u\|^{2}_{L^{2}(4 I\times 2Q)},$$ moreover, by Young's inequality, $$U_{52}\leq C(\ez) \|u\|^{2}_{L^{2}(4 I\times 2Q)}+\ez \|\nabla u\|^{2}_{L^{2}(4 I\times 2Q)}.$$ In summary, upon recalling \eqref{eq: 2.9}, we attain $$\|\nabla u\|^{2}_{L^{2}(I\times Q)}\leq C(\ez) \|u\|^{2}_{L^{2}(4 I\times 2Q)}+\ez \|\nabla u\|^{2}_{L^{2}(4 I\times 2Q)},$$ which yields \eqref{eq: 2.6} after scaling. 

Observing that $u\in V=L^{2}(\rr; W^{1, 2}(\rz)),$ by passing to the limit $r(Q)\to \infty$ in \eqref{eq: 2.6}, we derive $$ \|\nabla u\|^{2}_{L^{2}(\rdm)}\leq \ez\|\nabla u\|^{2}_{L^{2}(\rdm)}.$$ As a consequence, by letting $\ez\to 0,$  we get $\nabla u=0,$ which implies $u(t, x)=u(t).$ Furthermore, the fact that $u\in L^{2}(\rr^{n+1})$ leads to $u\equiv 0$ and $v=\tilde{v}\in E.$ From this, Lemma \ref{lemma: 2.3} and Lemma \ref{lemma: 2.42} deliver $$\|D_{t}^{1/2}v\|_{L^{2}(\rdm)}\lesssim \|\H_{\kappa} v\|_{E^{*}}\lesssim \|v\|_{L^{2}(\rdm)}+\|f\|_{L^{2_{*}}(\rdm)}+\|F\|_{L^{2}(\rdm)}.$$ Through a density argument and an application of Lemma \ref{lemma: 2.42}, it follows that  $v\in L^{2^{*}}(\rdm).$ Since $\partial_{t}v=\mbox{div}((S(t, x)+D(t, x))\nabla v)+f+\mbox{div} F,$ we finally arrive at \begin{equation}\label{eq: 2.10}\|\partial_{t}v\|_{V_{2^{*}}^{*}}\lesssim \|v\|_{V}+\|v\|_{L^{2^{*}}(\rdm)}+\|f\|_{L^{2_{*}}(\rdm)}+\|F\|_{L^{2}(\rdm)}.\end{equation}
With \eqref{eq: 2.10} in hand, Lemma \ref{lemma: 2.41} applies and the conclusion $(iii)$ follows. This completes the proof of Proposition \ref{pro: 2.4}.

\hfill$\Box$ 
	
\begin{remark}\label{rem: 2.5}\; If, in addition, $A\in C^{\infty}(\rdm),$ the estimate \eqref{eq: 2.6} can be strengthened to the classical Caccioppoli inequality: \begin{equation}\label{eq: 2.11}
\fint_{I\times Q}|\nabla u|^{2}\lesssim \frac{1}{r(Q)^{2}}\fint_{4I\times 2Q}| u|^{2}.\end{equation}
\end{remark} In fact, in this situation, we are able to test the equation \eqref{eq: 2.5} with $\eta^{2}u,$ where, for any $0<\rho<R<\infty,$ $\eta: \rr^{n+1}\to \rr$ is supposed to be a smooth function with $0\leq\rho\leq 1,$ $\eta\equiv 1$ on $I_{\rho^{2}}\times Q_{\rho},$ $\supp \eta \subset I_{\frac{(\rho+R)^{2}}{4}}\times Q_{\frac{\rho+R}{2}}$ and 
\begin{equation*}
|\nabla \eta|^{2}+|\nabla^{2} \eta|+|\partial_{t}\eta|\lesssim \frac{1}{(R-\rho)^{2}}.\end{equation*}
Using a standard argument in the theory of PDEs, one can easily deduce that
\begin{equation} \label{eq: 2.17}\int_{I_{\rho^{2}}\times Q_{\rho}}|\nabla u|^{2}\leq \frac{C}{(R-\rho)^{2}}\int_{ I_{R^{2}}\times Q_{R}}|u|^{2}+\frac{1}{2}\int_{ I_{R^{2}}\times Q_{R}}|\nabla u|^{2}\end{equation} holds for any $0<\rho<R<\infty.$ It should be particularly noted that, in the process of deriving \eqref{eq: 2.17}, the method developed in Proposition \ref{pro: 2.4} to settle $U_{5}$ must be borrowed once more to resolve the terms involving the anti-symmetric matrix $D.$ Once \eqref{eq: 2.17} is verified, \eqref{eq: 2.11} follows at once by employing the well-known iteration lemma, that is, 
 
\begin{lemma}\label{lemma: 2.6}(\cite[Lemma 5.1]{G})\; Let $f(t)$ be a nonnegative bounded function on $0\leq T_{0}\leq t\leq T_{1}.$ Suppose that for $ T_{0}\leq t<s\leq T_{1}$ we have $$f(t)\leq A(s-t)^{-\alpha}+B+\theta f(s)$$ with $\alpha>0,$ $0<\theta<1$ and $A, B$ nonnegative constants. Then there exists a constant $c=c(\alpha, \theta)$ such that for all $T_{0}\leq \rho<R\leq T_{1}$ we have $$f(\rho)\leq c[A(R-\rho)^{-\alpha}+B].$$
 \end{lemma} 
 
\begin{remark}\label{rem: 2.7}\; It is not clear whether the parabolic equation \eqref{eq: 2.5}, featuring partially unbounded coefficients and satisfying the G$\mathring{a}$rding inequality, can be approximated by a sequence of equations with smooth coefficients. In the proof of \cite[Theorem 1.1]{SS}, such an approximation appears to rely a priori on a (local) higher integrability (greater than 2) of the gradient. If this conclusion is valid, the restriction $n\geq 3$ imposed in Proposition \ref{pro: 2.4} can be removed by Remark \ref{rem: 2.5}. 
\end{remark} 
 
\section{Local regularity estimates for weak solutions in $L^{2}$}

This section is devoted to establishing the improved Caccioppoli inequality for weak solutions to parabolic systems considered here. In contrast to the case with $L^{\infty}$ coefficients, we are unable to prove a Caccioppoli inequality similar to \cite[Proposition 4.3]{ABES} at the outset. Instead, we obtain a weaker version (see Proposition \ref{pro: 3.2} below), fortunately, it is sufficient to achieve our goal. 

\begin{theorem}\label{thm: 3.1}\;  Set $\Omega:=I_{0}\times Q_{0}.$ Assume $f\in L_{loc}^{2_{*}}(\Omega, \cc^{m})$ and $F\in L^{2}_{loc}(\Omega, \cc^{m}),$ if $u$ is a weak solution of $\partial_{t}u-\mbox{div}((S(t, x)+D(t, x))\nabla u)=f+\mbox{div} F$ in $\Omega,$ then $$u\in  L_{loc}^{2^{*}}(\Omega, \cc^{m})\cap C_{0}(I_{0}; L^{2}_{loc}(Q_{0})).$$ More precisely, for any $\eta\in C_{0}^{\infty}(\Omega),$ the function $t\to \|(u\eta)(t, \cdot)\|_{L^{2}(\rz)}^{2}$ is absolutely continuous on $\rr,$ and $D_{t}^{1/2}(u\eta)\in L^{2}(\rr^{n+1}; \cc^{m}),$ in particular, \begin{equation}\label{eq: 3.1}\|u\eta\|_{L^{2^{*}}(\rdm)}\lesssim \|u\eta\|_{L^{2}(\rdm)}+\|\tilde{f}\|_{L^{2_{*}}(\rdm)}+\|\tilde{F}\|_{L^{2}(\rdm)}+\|\nabla(u\eta)\|_{L^{2}(\rdm)},\end{equation} where \begin{equation}\label{eq: 3.2}\begin{aligned}\tilde{f}&:=[u\eta+(\partial_{t}\eta)u-S\nabla u\cdot \nabla \eta-F\cdot \nabla \eta]-\l(D-\fint_{Q_{0}}D\r)\nabla u\cdot \nabla \eta\\
&:=\tilde{f}_{S}-\l(D-\fint_{Q_{0}}D\r)\nabla u\cdot \nabla \eta\end{aligned}\end{equation} and \begin{equation}\label{eq: 3.3}\tilde{F}:=[-S(u\nabla \eta)+F\eta]-\l(D-\fint_{Q_{0}}D\r)(u\nabla \eta):=\tilde{F}_{S}-\l(D-\fint_{Q_{0}}D\r)(u\nabla \eta).\end{equation}
The implicit constant depends only on dimensions, $\|D\|_{L^{\infty}(\rr; BMO(\rz))}$ and ellipticity. 
\end{theorem}

{\it Proof.}\quad Given any $\eta \in C_{0}^{\infty}(\Omega).$ Note that $u$ is in priori in $W^{1, 2}_{loc}(\Omega)$ by the definition of weak solutions. Then, in view of \cite[Section 3]{XQ} or \cite[Remark 2.8]{AK}, we see that $u$ is a weak solution of $\partial_{t}u-\mbox{div}\l(\l(S(t, x)+\l(D-\fint_{Q_{0}}D\r)\r)\nabla u\r)=f+\mbox{div} F$ in $\Omega.$ From this, a straightforward calculation reveals that $v:=u\eta\in V$ is a weak solution of \begin{equation}\label{eq: 3.4}\partial_{t}v-\mbox{div}\l(\l(S(t, x)+\l(D-\fint_{Q_{0}}D\r)\r)\nabla v\r)=\tilde{f}+\mbox{div} \tilde{F}\quad\mbox{in}\;\rdm\end{equation} with $\tilde{f}$ and $\tilde{F}$ defined by \eqref{eq: 3.2} and \eqref{eq: 3.3}, respectively. It is evident that $\l(D-\fint_{Q_{0}}D\r)$ retains the properties of $D$ and \begin{equation}\label{eq: 3.5}\l(D-\fint_{Q_{0}}D\r)\in L^{q}_{loc}(\rdm) \quad \mbox{for any }\; q>1.\end{equation} Using H\"{o}lder inequality ($2_{*}<2$) in combination with \eqref{eq: 3.5}, one can readily demonstrate that $\tilde{f}\in L^{2_{*}}(\rdm)$ and $\tilde{F}\in L^{2}(\rdm).$ The proof follows a similar line of reasoning as Proposition \ref{pro: 3.2} below, and we leave the details to the reader. Therefore, \eqref{eq: 3.1} is an immediate consequence of Proposition \ref{pro: 2.4}.

\hfill$\Box$ 

Proposition \ref{pro: 3.2} is derived from the combination of Theorem \ref{thm: 3.1} and Proposition \ref{pro: 2.4}. 

\begin{proposition}\label{pro: 3.2}\; Let $f\in L_{loc}^{2_{*}}(\Omega, \cc^{m})$, $F\in L^{2}_{loc}(\Omega, \cc^{m})$ and $B\in L^{\infty}(\Omega; \L(\cc^{m}))$ satisfying \begin{equation}\label{eq: 3.6}\mbox{Re}\; \int_{\Omega}B\varphi\cdot \overline{\varphi}\geq 0\quad \mbox{for any}\; \;\varphi\in L^{2}(\Omega; \cc^{m}).\end{equation}  Assume that $u$ is a weak solution of $\partial_{t}u-\mbox{div}((S(t, x)+D(t, x))\nabla u)=f+\mbox{div} F-Bu$ in $\Omega.$  Let $I\times Q$ be an open parabolic cylinder with $l(I)\approx r(Q)^{2}$ such that $\overline{\gamma^{2}I\times \gamma Q}\subset \Omega$ for some $\gamma>1.$ Thus for any $s>2$ \begin{equation}\label{eq: 3.7}\begin{aligned}\l(\fint_{I\times Q}|\nabla u|^{2}\r)^{1/2}&\lesssim \frac{1}{r(Q)}\l(\fint_{\gamma^{2}I\times \gamma Q}|u|^{s}\r)^{1/s}+\l(\fint_{\gamma^{2}I\times \gamma Q}|F|^{2}\r)^{1/2}\\
&\quad\quad+r(Q)\l(\fint_{\gamma^{2}I\times \gamma Q}|f|^{2_{*}}\r)^{1/2_{*}}.\end{aligned}\end{equation} 
The implicit constant depends only on dimensions, $\|D\|_{L^{\infty}(\rr; BMO(\rz))}$, $\gamma,$ ellipticity, the constants controlling the ratio $l(I)/r(Q)^{2}$ and $\|B\|_{\infty}.$
\end{proposition}

\begin{remark}\label{rem: 3.3}\; Although the exponent $s$ in \eqref{eq: 3.7} may be arbitrarily close to $2,$ the proof presented here does not accommodate the case $s=2$. However, it is eventually proven that $s$ can be any value in the range $[1, 2]$. \end{remark}

{\it Proof.}\quad In order to simplify the exposition, we assume that $r(Q)=1$ and $\gamma=2$. Select $\eta\in C_{0}^{\infty}(4I\times 2Q)$ such that  $0\leq \eta\leq 1$ and $\eta\equiv 1$ on $I\times Q,$ and let $v:=u\eta.$ Theorem \ref{thm: 3.1} directly yields that $v$ is a weak solution of $\partial_{t}v-\mbox{div}((S(t, x)+\l(D-\fint_{2Q}D\r))\nabla v)=\tilde{f}+\mbox{div} \tilde{F}-Bv\quad\mbox{in}\;\rdm,$ where $\tilde{f}, \tilde{F}$ are respectively defined by \eqref{eq: 3.2} and \eqref{eq: 3.3} with $2Q$ in lieu of $Q_{0}.$ Owing to the assumptions on $f, F, B,$ we can easily deduce that $\tilde{f}-Bv\in L^{2_{*}}(\rdm)$ and $\tilde{F}\in L^{2}(\rdm).$ Thus $v\in V_{2^{*}}$ thanks to Proposition \ref{pro: 2.4}. Furthermore, the same reasoning applied in the proof of \eqref{eq: 2.10} allows us to conclude that $\partial_{t}v\in V_{2^{*}}^{*}.$ From the above analysis, we have
\begin{equation*}\begin{aligned}
 0=&-\mbox{Re}\;\int_{\rdm}S\nabla v\cdot \overline{\nabla v}-\tilde{F}_{S}\cdot \overline{\nabla v}-Re\;\int_{\rdm}\l(D-\fint_{2Q}D\r)\nabla v \cdot \overline{\nabla v}\\
&\quad +\mbox{Re}\;\int_{\rdm}\l(D-\fint_{2Q}D\r)(u\nabla \eta)\cdot \overline{\nabla v}+Re\;\int_{\rdm}\tilde{f}_{S}\overline{ v}\\
&\quad\quad-\mbox{Re}\;\int_{\rdm}\l(D-\fint_{2Q}D\r)\nabla u\cdot \nabla \eta\overline{ v}\\
&\quad\quad-\mbox{Re}\;\int_{\rdm} Bv\cdot \overline{ v}:=-J_{1}-J_{2}-J_{3}+J_{4}+J_{5}-J_{6}-J_{7}
\end{aligned} \end{equation*} taking into account of Lemma \ref{lemma: 2.41}. First, by \eqref{eq: 2.3}, $$J_{3}=0.$$ Second, using \eqref{eq: 3.6}, $$-J_{7}\leq 0,$$ hence \begin{equation}\label{eq: 3.8}\mbox{Re}\;\int_{\rdm}S\nabla v\cdot \overline{\nabla v}=J_{1}\leq |J_{2}|+|J_{4}|+|J_{5}|+|J_{6}|.\end{equation} Concatenating \eqref{eq: 1.2} and \eqref{eq: 3.8}, we gain \begin{equation}\label{eq: 3.9}\int_{\rdm}|\nabla v|^{2}\lesssim  \int_{\rdm}|u|^{2}+|J_{2}|+|J_{4}|+|J_{5}|+|J_{6}|.\end{equation}

We proceed by eliminating $J_{6}$ in the sequel. \begin{equation*}\begin{aligned}
 J_{6}&=\mbox{Re}\;\sum_{\alpha=\beta}\int_{\rdm}\l(D_{kj}^{\alpha \beta}-\fint_{2Q}D_{kj}^{\alpha \beta}\r)\partial_{k}u^{\alpha}\partial_{j}\eta\eta\overline{u^{\beta}}\\
 &=\mbox{Re}\;\sum_{\alpha=\beta}\int_{\rdm}\l(D_{kj}^{\alpha \beta}-\fint_{2Q}D_{kj}^{\alpha \beta}\r)\partial_{k}(u^{\alpha}\eta)\partial_{j}\eta\overline{u^{\beta}}\\
 &\quad\quad-\mbox{Re}\;\sum_{\alpha=\beta}\int_{\rdm}\l(D_{kj}^{\alpha \beta}-\fint_{2Q}D_{kj}^{\alpha \beta}\r)\partial_{k}\eta\partial_{j}\eta\overline{u^{\beta}}u^{\alpha}
   \end{aligned}\end{equation*}
 \begin{equation*}\begin{aligned}  
&=\mbox{Re}\;\sum_{\alpha=\beta}\int_{\rdm}\l(D_{kj}^{\alpha \beta}-\fint_{2Q}D_{kj}^{\alpha \beta}\r)\partial_{k}(u^{\alpha}\eta)\partial_{j}\eta\overline{u^{\beta}}\end{aligned}\end{equation*} where we have used the fact that \begin{equation}\label{eq: 3.10}\mbox{Re}\;\sum_{\alpha=\beta}\int_{\rdm}\l(D_{kj}^{\alpha \beta}-\fint_{2Q}D_{kj}^{\alpha \beta}\r)\partial_{k}\eta\partial_{j}\eta\overline{u^{\beta}}u^{\alpha}=0.\end{equation}
Since \eqref{eq: 3.5} and $s>2,$ there exists a $q_{s}>2$ such that $\frac{1}{q_{s}}+\frac{1}{2}+\frac{1}{s}=1.$ Utilizing H\"{o}lder inequality and Young's inequality, we achieve \begin{equation*}\begin{aligned}
 |J_{6}|&\lesssim \bigg|\mbox{Re}\;\sum_{\alpha=\beta}\int_{\rdm}\l(D_{kj}^{\alpha \beta}-\fint_{2Q}D_{kj}^{\alpha \beta}\r)\partial_{k}(u^{\alpha}\eta)\partial_{j}\eta\overline{u^{\beta}}\bigg|\\
&\lesssim \l(\int_{4I}\int_{2Q}\bigg|D(t, \cdot)-\fint_{2Q}D\bigg|^{q_{s}}\r)^{1/q_{s}}\l(\int_{4I\times 2Q}|\nabla v|^{2}\r)^{1/2}\l(\int_{4I\times 2Q}|u\nabla \eta|^{s}\r)^{1/s}\\
&\lesssim \|\|D(t, \cdot)\|_{BMO(\rz)}\|_{L^{\infty}(\rr)}\l(\int_{\rdm}|\nabla v|^{2}\r)^{1/2}\l(\int_{\rdm}|u\nabla \eta|^{s}\r)^{1/s}\\
&\lesssim \ez \int_{\rdm}|\nabla v|^{2}+\l(\int_{4I\times 2Q}|u|^{s}\r)^{2/s}.
\end{aligned}\end{equation*}  The term $J_{4}$ can be handled by the same approach, hence we omit the details. By invoking Young's inequality again and letting $\ez$ small enough, we get \begin{equation}\label{eq: 3.11}\begin{aligned}\int_{\rdm}|\nabla v|^{2}&\lesssim \int_{\rdm}|\tilde{F}_{S}|^{2}+|J_{5}|+\ez \int_{\rdm}|\nabla v|^{2}+c(\ez)\l(\int_{\rdm}|u\nabla \eta|^{s}\r)^{2/s}\\
&\lesssim \int_{\rdm}|\tilde{F}_{S}|^{2}+|J_{5}|+\l(\int_{4I\times 2Q}|u|^{s}\r)^{2/s}\end{aligned}\end{equation} Next,
we use the methods from \cite{ABES} to analyze the first two terms on the right-hand side of \eqref{eq: 3.11}. Compared to parabolic systems with bounded coefficients, several technical details need to be emphasized during the proof.

Pick another real-valued function $\psi\in C_{0}^{\infty}(4I\times 2Q)$ such that $\psi\equiv 1$ on $\supp \eta.$ Then, the ingredient in $J_{5}$ can be rewritten as $$\tilde{f}_{S}\overline{v}=f\psi\eta^{2}\overline{u}+\partial_{t}\eta\eta|u|^{2}-S\nabla u\cdot \nabla\eta \eta \overline{u}-F\cdot\nabla \eta\eta\overline{u}.$$ By recalling the definitions of $\tilde{F}_{S}$ ($S\in L^{\infty}(\rdm)$) and $J_{5}$, and inserting the above equality into \eqref{eq: 3.11}, it is straightforward to verify that \begin{equation*}\begin{aligned}
\int_{\rdm}|\nabla v|^{2}&\lesssim L+\int_{\rdm}|f\psi\eta^{2}\overline{u}|+|S\nabla u\cdot \nabla\eta \eta \overline{u}|+\l(\int_{\rdm}|u\nabla \eta|^{s}\r)^{2/s}\\
&:=L+L_{1}+L_{2}+\l(\int_{4I\times 2Q}|u|^{s}\r)^{2/s},
\end{aligned}\end{equation*} where $$L:=\int_{\rdm}|u|^{2}(|\nabla \eta|^{2}+|\partial_{t}\eta|+\eta^{2}).$$  Since $S\in L^{\infty},$ we can apply Young's inequality to estimate $L_{2}$ as follows: $$L_{2}\lesssim\ez\int_{\rdm}|\nabla v|^{2}+L.$$ Consequently, \begin{equation}\label{eq: 3.12}\int_{\rdm}|\nabla v|^{2}\lesssim L+\l(\int_{4I\times 2Q}|u|^{s}\r)^{2/s}+L_{1}.\end{equation} It remains to bound $L_{1}.$ To the end, by setting $\omega:=u\eta^{2},$ we conclude from Theorem \ref{thm: 3.1} that $\omega$ is a weak solution of $$\partial_{t}\omega-\mbox{div}\l(\l(S(t, x)+\l(D-\fint_{2Q}D\r)\r)\nabla \omega\r)=\tilde{f}'+\mbox{div} \tilde{F}'-B\omega \quad\mbox{in}\;\rdm,$$ where $$\tilde{f}':=\tilde{f}_{S}^{\eta^{2}}-\l(D-\fint_{2Q}D\r)\nabla u\cdot\nabla\eta^{2}, \quad\tilde{F}':=\tilde{F}_{S}^{\eta^{2}}-\l(D-\fint_{2Q}D\r)u\nabla \eta^{2},$$ and $\tilde{f}_{S}^{\eta^{2}}, \tilde{F}_{S}^{\eta^{2}}$ are defined in the same way as $\tilde{f}_{S}, \tilde{F}_{S}$ in Theorem \ref{thm: 3.1}, but with $\eta$ substituted by $\eta^{2}.$ Applying \eqref{eq: 2.4} in Proposition \ref{pro: 2.4} and crudely using H\"{o}lder inequality($2_{*}<2$), we obtain 
\begin{align*}
\|\omega\|_{L^{2^{*}}(\rdm)}&\lesssim \|\tilde{f}'-B\omega\|_{L^{2_{*}}(\rdm)}+\|\tilde{F}'\|_{L^{2}(\rdm)}+\|\omega\|_{V}\\[4pt]
&\lesssim \|\tilde{f}'\|_{L^{2_{*}}(\rdm)}+\|\tilde{F}'\|_{L^{2}(\rdm)}+\|\omega\|_{L^{2}(\rdm)}+\|\nabla \omega\|_{L^{2}(\rdm)}.\end{align*}  
Note that \begin{align*} \|\l(D-\fint_{2Q}D\r)u\nabla\eta^{2}\|_{L^{2}(\rdm)}&\lesssim\l(\int_{4I\times 2Q}\bigg|\l(D-\fint_{2Q}D\r)\bigg|^{2\l(\frac{s}{2}\r)'}\r)^{\frac{1}{2\l(\frac{s}{2}\r)'}}\|u\nabla \eta^{2}\|_{L^{s}}\\[4pt]
&\lesssim \|\|D(t, \cdot)\|_{BMO(\rz)}\|_{L^{\infty}(\rr)}\|u\|_{L^{s}(4I\times 2Q)}.\end{align*} 
By the same token, 
\begin{align*}
 \|\l(D-\fint_{2Q}D\r)\nabla u\cdot\nabla\eta^{2}\|_{L^{2_{*}}(\rdm)}&\lesssim \|\l(D-\fint_{2Q}D\r)\nabla v\cdot\nabla\eta\|_{L^{2_{*}}(\rdm)}\\[4pt]
 &\quad + \|\l(D-\fint_{2Q}D\r)u\nabla \eta\cdot\nabla\eta\|_{L^{2_{*}}(\rdm)}\\[4pt]
 &\lesssim \l(\int_{4I\times 2Q}\bigg|\l(D-\fint_{2Q}D\r)\bigg|^{2_{*}\l(\frac{2}{2_{*}}\r)'}\r)^{\frac{1}{2_{*}\l(\frac{2}{2_{*}}\r)'}}\|\nabla v\|_{L^{2}}\\[4pt]
 &\quad\quad+\|\|D(t, \cdot)\|_{BMO(\rz)}\|_{L^{\infty}(\rr)}\|u\|_{L^{s}(4I\times 2Q)}\\[4pt]
 &\lesssim \|\|D(t, \cdot)\|_{BMO(\rz)}\|_{L^{\infty}(\rr)}(\|\nabla v\|_{L^{2}}+\|u\|_{L^{s}(4I\times 2Q)}).\end{align*}  

In conclusion, we get 
\begin{equation*}\begin{aligned}
&\int_{\rdm}|f\psi\eta^{2}\overline{u}|\lesssim \|f\psi\|_{L^{2_{*}}(\rdm)}\|\omega\|_{L^{2^{*}}(\rdm)}\\
&\quad\lesssim \|f\psi\|_{L^{2_{*}}(\rdm)}\\
&\quad \times\l(\|\tilde{f}_{S}^{\eta^{2}}\|_{L^{2_{*}}(\rdm)}+\|u\|_{L^{s}(4I\times 2Q)}+\|\tilde{F}_{S}^{\eta^{2}}\|_{L^{2}(\rdm)}+\|\nabla v\|_{L^{2}}+\|\omega\|_{L^{2}(\rdm)}+\|\nabla \omega\|_{L^{2}(\rdm)}\r)\\
&\lesssim L+\|f\psi\|_{L^{2_{*}}(\rdm)}^{2}+\|f\psi\|_{L^{2_{*}}(\rdm)}\|S\nabla u\cdot \nabla \eta^{2}\|_{L^{2}}\\
&\quad\quad+\ez\|\nabla v\|_{L^{2}}^{2}+\|u\|_{L^{s}(4I\times 2Q)}^{2}+\|f\psi\|_{L^{2_{*}}(\rdm)}\|\nabla \omega\|_{L^{2}},
\end{aligned}\end{equation*} taking also into account that $$\|\tilde{f}_{S}^{\eta^{2}}\|_{L^{2_{*}}(\rdm)}\lesssim \|f\eta\|_{L^{2_{*}}(\rdm)}+\|\partial_{t}\eta^{2}u\|_{L^{2}(\rdm)}+\|S\nabla u\cdot \nabla \eta^{2}\|_{L^{2}}+\|F\nabla \eta^{2}\|_{L^{2}(\rdm)}$$ and  $$\|\tilde{F}_{S}^{\eta^{2}}\|_{L^{2}(\rdm)} \lesssim \|S(u\nabla \eta^{2})\|_{L^{2}(\rdm)}+\|F\eta^{2}\|_{L^{2}(\rdm)}.$$ For the estimates of $\|S\nabla u\cdot \nabla \eta^{2}\|_{L^{2}}$ and $\|\nabla \omega\|_{L^{2}},$ we may run a similar argument as in the proof of \cite[Proposition 4.3]{ABES}. To save pages, we exclude the details. Ultimately, we reach the conclusion that $$L_{1}\lesssim L+\|f\psi\|_{L^{2_{*}}(\rdm)}^{2}+\|u\|^{2}_{L^{s}(4I\times 2Q)}.$$ This confirms \eqref{eq: 3.7}, thereby completing the proof. 

\hfill$\Box$ 

Theorem \ref{thm: 3.1} and Proposition \ref{pro: 3.2} together imply the following reverse H\"{o}lder inequality for solutions to \eqref{eq: 1.1}. 

\begin{proposition}\label{pro: 3.4}\; Let $f\in L_{loc}^{2_{*}}(\Omega, \cc^{m})$ and $F\in L^{2}_{loc}(\Omega, \cc^{m}).$ If $u$ is a weak solution of $\partial_{t}u-\mbox{div}((S(t, x)+D(t, x))\nabla u)=f+\mbox{div} F$ in $\Omega$, then we have \begin{equation}\label{eq: 3.13}\l(\fint_{I\times Q}| u|^{2^{*}}\r)^{1/2^{*}}\lesssim \fint_{\gamma^{2}I\times \gamma Q}|u|+r(Q)\l(\fint_{\gamma^{2}I\times \gamma Q}|F|^{2}\r)^{1/2}+r(Q)^{2}\l(\fint_{\gamma^{2}I\times \gamma Q}|f|^{2_{*}}\r)^{1/2_{*}}.\end{equation} 
The implicit constant depends only on dimensions, $\|D\|_{L^{\infty}(\rr; BMO(\rz))}$, $\gamma,$ ellipticity, and the constants controlling the ratio $l(I)/r(Q)^{2}$.
\end{proposition}

{\it Proof.}\quad As before, we assume $r(Q)=1.$ Then, by \eqref{eq: 3.1} and \eqref{eq: 3.7} \begin{equation}\label{eq: 3.14}\begin{aligned}
\l(\fint_{I\times Q}| u|^{2^{*}}\r)^{1/2^{*}}&\lesssim \l(\int_{\rdm}| u\eta|^{2^{*}}\r)^{1/2^{*}}\\
&\lesssim \|u\eta\|_{L^{2}(\rdm)}+\|\tilde{f}\|_{L^{2_{*}}(\rdm)}+\|\tilde{F}\|_{L^{2}(\rdm)}+\|\nabla(u\eta)\|_{L^{2}(\rdm)}\\
&\lesssim \l(\fint_{\frac{9}{4}I\times \frac{3}{2} Q}|u|^{s}\r)^{1/s}+\l(\fint_{\frac{9}{4}I\times \frac{3}{2} Q}|F|^{2}\r)^{1/2}+\l(\fint_{\frac{9}{4}I\times \frac{3}{2} Q}|f|^{2_{*}}\r)^{1/2_{*}}\\
\end{aligned} \end{equation} 
\begin{equation*}\begin{aligned}
&\quad\quad\quad\quad\quad\quad\quad\quad\quad\quad\quad+\l(\fint_{\frac{9}{4}I\times \frac{3}{2} Q}|\nabla u|^{2}\r)^{1/2}\\
&\quad\quad\quad\quad\quad\quad\quad\quad\quad\lesssim \l(\fint_{4I\times 2 Q}|u|^{s}\r)^{1/s}+\l(\fint_{4I\times 2Q}|F|^{2}\r)^{1/2}+\l(\fint_{4I\times 2Q}|f|^{2_{*}}\r)^{1/2_{*}},
\end{aligned} \end{equation*} where in the last second step we have adopted a similar argument as in the proof of Proposition \ref{pro: 3.2}. By \eqref{eq: 3.14}, if we pick $2<s<2^{*},$ \eqref{eq: 3.13} is a straightforward  consequence of the self-improved lemma \cite[Theorem B.1]{BCF}, that is, 

\begin{theorem}\label{thm: 3.5}\;  Let $(M, d, \mu)$ be a doubling metric measure space and $\mathcal{F}$ be the collection of all balls of the space $M.$ We consider a functional $a: \mathcal{F}\to [0, \infty)$ satisfying that there exists a constant $c > 0$ such that for every pair of balls $B, \tilde{B}$ with $B\subset  \tilde{B}\subset 2B$ $$ca(B)\leq a(\tilde{B})\leq c^{-1}a(2B).$$ Let $1<p<q<\infty,$ and we assume that $\omega\in L^{1}_{loc}(M, \mu)$ is a non-negative function such that for every $B\subset M$ $$\l(\fint_{B}\omega^{q}d\mu\r)^{1/q}\lesssim\l(\fint_{2B}\omega^{p}d\mu\r)^{1/p}+a(B),$$ then for any $\eta\in (0, 1)$ and every ball $B\subset M$ $$\l(\fint_{B}\omega^{q}d\mu\r)^{1/q}\lesssim\l(\fint_{2B}\omega^{\eta p}d\mu\r)^{1/(\eta p)}+a(2B).$$
\end{theorem}
\hfill$\Box$ 

From the previous two propositions,  we construct the improved Caccioppoli inequality, as stated in Corollary \ref{cor: 3.6}, in agreement with the result in \cite{ABES}.

\begin{corollary}\label{cor: 3.6}\;  Under the same assumptions and notations as in Proposition \ref{pro: 3.2} and \ref{pro: 3.4}, we have $$\l(\fint_{I\times Q}|\nabla u|^{2}\r)^{1/2}\lesssim \frac{1}{r(Q)}\fint_{\gamma^{2}I\times \gamma Q}|u|+\l(\fint_{\gamma^{2}I\times \gamma Q}|F|^{2}\r)^{1/2}+r(Q)\l(\fint_{\gamma^{2}I\times \gamma Q}|f|^{2_{*}}\r)^{1/2_{*}}.$$ 
The implicit constant depends only on dimensions, $\|D\|_{L^{\infty}(\rr; BMO(\rz))},$ $\gamma,$ ellipticity, and the constants controlling the ratio $l(I)/r(Q)^{2}$.
\end{corollary}
 
\section{Proofs of Theorem \ref{thm: 4.1} and Theorem \ref{thm: 4.2}}

It is worth reiterating that the new insight in \cite[Theorem 8.1]{ABES} and \cite[Theorem 8.2]{ABES} lies in the weaker assumption $f\in L_{loc}^{2_{*}}(\Omega, \cc^{m})$ (where $2_{*}<2$), in contrast to earlier results. This novelty is maintained in the present scenario. 

\begin{proof}[Proofs of Theorem \ref{thm: 4.1} and Theorem \ref{thm: 4.2}]
\quad We particularly stress that all the principal results in \cite[Section 6; Section 7; Section 8]{ABES} have parallel versions for parabolic systems having a skew-symmetric and real BMO part, where the implicit constants depend unavoidably on dimensions, $\|D\|_{L^{\infty}(\rr;  BMO(\rz))}$ and ellipticity. Furthermore, their proofs can essentially be carried over from \cite[Section 6; Section 7; Section 8]{ABES}, with Proposition \ref{pro: 2.4} serving as a key component. Although we occasionally rely on the equation \eqref{eq: 1.1} in the proof, indicating that the anti-symmetric matrix $D$ has an influence, the additional terms involving $D$-compared to the case with $L^{\infty}$ coefficients-can be easily managed by using the techniques developed in the preceding sections. Simply recall to use Proposition \ref{pro: 2.1} and Proposition \ref{pro: 3.4} when necessary. Working out the details along the process of proof in \cite{ABES} is left to interested readers. \end{proof}

Theorem \ref{thm: 4.2} represents the reverse H\"{o}lder inequality for the gradient in a non-homogeneous context. It also generalizes \cite[Theorem 2.2]{SS} which dealt with parabolic equations with $f=F=0$ under the assumption that $S\in L^{\infty}$ is real and uniformly elliptic and $D\in L^{\infty}(BMO)$ is real and skew-symmetric. For completeness, we present here an alternative approach, as mentioned in \cite[Remark 8.3]{ABES}, to prove the result. We highlight that this approach is rooted in the classical reasoning of \cite[Lemma 1.2]{GS} or \cite[Chapter \;V \;Theorem\; 2.2]{GM}, and offers a more direct proof compared to that in \cite{ABES}.

\begin{proof}[An alternative proof for Theorem \ref{thm: 4.2}]\quad  We first consider $16I\times 4Q\subset \Omega$ with $I:=(-1, 0)$ and $Q:=B(0, 1).$ Then $4I:=(-4, 0)$ and $2Q:=B(0, 2).$ Our purpose in the subsequent discussion is to build \begin{equation}\label{eq: 4.2}\begin{aligned}
\int_{I\times Q}|\nabla v|^{2}&\lesssim\|f\|_{L^{2_{*}}( 16I\times 4Q)}^{2}+\l(\int_{16I\times 4Q}|F|^{2}\r)^{2/2}+\l(\int_{16I\times 4Q}|\nabla v|^{r}\r)^{2/r}\\
&\quad\quad+\ez \int_{16I\times 4Q}|\nabla v|^{2}.\end{aligned}\end{equation}  

Prior to starting, we make some preliminary preparations. First, we let $\eta:=\eta_{I}(t)\eta_{Q}(x)$ such that $\eta_{Q}\in C_{0}^{\infty}(2Q),$ $\eta_{Q}\equiv 1$ on $Q$ and $\eta_{I}(t)\equiv 1$ when $t\geq -1,$ in particular, $\eta_{I}(t)=0$ when $t\leq - 4.$ Second, we set $$\hat{v}(t, x):=v(t, x)-\tilde{v}(t):=v(t, x)-\frac{\int_{\rz}v(t, x)\eta_{Q}(x)^{2}dx}{\int_{\rz}\eta_{Q}(x)^{2}}.$$
 
Armed with the estimate $(iii)$ (especially the expression for $D_{t}^{1/2}$) in Proposition \ref{pro: 2.4}, we are able to test the equation with $\varphi:=\hat{v}(t, x)\eta^{2}1_{(-\infty, t_{0})}$ for any $-2<t_{0}<0$ (see \cite{GS}). Through integration by parts, we can deduce the following inequality: \begin{equation*}\begin{aligned}
\frac{1}{2}\int_{2Q}|\hat{v}(t_{0}, x)|^{2}\eta_{Q}^{2}&+\mbox{Re}\;\int_{\rdm\cap (-\infty, t_{0})}\eta_{I}^{2}S\nabla(\eta_{Q}\hat{v})\overline{\nabla(\eta_{Q}\hat{v})}\\
&\quad\quad\leq \mbox{Re}\;\int_{\rdm\cap (-\infty, t_{0})}\eta^{2}f\cdot \overline{\hat{v}}-F\cdot \overline{\nabla(\hat{v}\eta^{2})}-S\nabla \hat{v}\cdot \eta_{Q}\nabla \eta_{Q}\overline{\hat{v}}\eta_{I}^{2}\\
&\quad\quad\quad-\mbox{Re}\;\int_{\rdm\cap (-\infty, t_{0})}\l(D-\fint_{2Q}D\r)\nabla \hat{v}\nabla \eta_{Q}^{2}\overline{\hat{v}}\eta_{I}^{2}+\frac{1}{2}|\hat{v}|^{2}\eta_{Q}^{2}\partial_{t}\eta_{I}^{2}\\
&\quad\quad\quad+ \mbox{Re}\;\int_{\rdm\cap (-\infty, t_{0})}\eta_{I}^{2}S\nabla \eta_{Q}\cdot \hat{v}\overline{\nabla\hat{v}}\eta_{Q}+\eta_{I}^{2}S\nabla\eta_{Q} \hat{v}\overline{\hat{v}}\nabla\eta_{Q}\\[4pt]
&\quad\quad:=
S_{1}+S_{2}+S_{3}+S_{4}+S_{5}+S_{6}+S_{7},\end{aligned}\end{equation*} where in the process we have also utilized two elementary facts, stated as  $$\int_{\rr\cap (-\infty, t_{0})}\l(\int_{\rz}\hat{v}\eta_{Q}^{2}\r)\partial_{t}\tilde{v}\eta_{I}^{2}=0$$ and $$\mbox{Re}\;\int_{\rdm\cap (-\infty, t_{0})}\l(\fint_{2Q}D\r)\nabla \hat{v}\cdot \nabla \varphi \overline{\hat{v}}=0, \quad \forall\; \varphi \in C_{0}^{\infty}(4I\times 2Q).$$ The two facts can be derived directly from the definition of $\hat{v}$ and the assumption on $D.$ 

We proceed by decomposing $S_{3}$ into \begin{equation*}\begin{aligned}S_{3}&=\mbox{Re}\;\int_{\rdm\cap (-\infty, t_{0})}S\nabla (\hat{v}\eta_{Q})\cdot \nabla \eta_{Q}\overline{\hat{v}}\eta_{I}^{2}-\mbox{Re}\;\int_{\rdm\cap (-\infty, t_{0})}S\nabla \eta_{Q}\hat{v}\cdot \nabla \eta_{Q}\overline{\hat{v}}\eta_{I}^{2}\\
&:=S_{31}+S_{32}.\end{aligned}\end{equation*} By directly applying Young's inequality and H\"{o}lder inequality, it can be proved that, for any $1<s<2,$ \begin{equation*}\begin{aligned}|S_{32}|&\lesssim \delta\int_{\rdm}|\nabla(\eta_{Q}\hat{v})|^{2}\eta_{I}^{2}+\int_{4I\times 2Q}|\hat{v}|^{2}\\
&\lesssim \delta\int_{\rdm}|\nabla(\eta_{Q}\hat{v})|^{2}\eta_{I}^{2}+\int_{4I}\l(\int_{2Q}|\hat{v}|^{\frac{2s}{2-s}}\r)^{\frac{2-s}{s}}\end{aligned}\end{equation*} and $$|S_{31}|+|S_{7}|\lesssim \int_{4I\times 2Q}|\hat{v}|^{2}\lesssim \int_{4I}\l(\int_{2Q}|\hat{v}|^{\frac{2s}{2-s}}\r)^{\frac{2-s}{s}}.$$ As the estimates for $S_{6}$ and $S_{32}$ are similar, and the same holds for $S_{5}$ and $S_{7},$ we omit them for brevity. Owing to $F\in L^{2}(4I\times 2Q),$ the term $S_{2}$ admits the following bound: $$|S_{2}|\lesssim \l(\int_{4I\times 2Q}|F|^{2}\r)^{2/2}+\delta \|\nabla (\hat{v}\eta)\|_{L^{2}(\rdm\cap (-\infty, t_{0}))}^{2}+\int_{4I}\l(\int_{2Q}|\hat{v}|^{\frac{2s}{2-s}}\r)^{\frac{2-s}{s}}.$$ We now turn to the analysis of $S_{4}.$ We split $S_{4}$ into \begin{equation*}\begin{aligned}S_{4}&=2\mbox{Re}\;\int_{\rdm}\l(D-\fint_{2Q}D\r)\nabla \l(\eta_{Q}\hat{v}\r)\cdot \nabla\eta_{Q}\overline{\hat{v}}\eta_{I}^{2}-2\mbox{Re}\;\int_{\rdm}\l(D-\fint_{2Q}D\r)\nabla \eta_{Q}\hat{v}\cdot \nabla\eta_{Q}\overline{\hat{v}}\eta_{I}^{2}\\
&:=S_{41}+S_{42}. \end{aligned}\end{equation*} The proof of $S_{41}$ follows directly from the argument presented at the beginning of \cite[Page 516]{SS}. Since there is no significant difference, we state the result without further elaboration, that is $$|S_{41}|\lesssim \delta \|\nabla (\hat{v}\eta)\|_{L^{2}(\rdm\cap (-\infty, t_{0}))}^{2}+\int_{4I}\l(\int_{2Q}|\hat{v}|^{\frac{2s}{2-s}}\r)^{\frac{2-s}{s}}.$$
Moreover, for any $1<s<2,$ \begin{equation*}\begin{aligned}
|S_{42}|&\lesssim \int_{-4}^{0}\eta_{I}^{2}\l(\int_{2Q}\bigg|D-\fint_{2Q}D\bigg|^{\frac{s}{2s-2}}\r)^{\frac{2s-2}{s}}\l(\int_{2Q}|\hat{v}|^{\frac{2s}{2-s}}\r)^{\frac{2-s}{s}}\\
&\lesssim \|\|D(t, \cdot)\|_{BMO(\rz)}\|_{L^{\infty}(\rr)}\int_{-4}^{0}\l(\int_{2Q}|\hat{v}|^{\frac{2s}{2-s}}\r)^{\frac{2-s}{s}}.
\end{aligned}\end{equation*}
It remains to estimate $S_{1}$ which is quite important as it captures the main distinctions from the homogeneous case in \cite{SS}. Unlike $S_{2}, S_{3}, ..., S_{7},$ the term $S_{1}$ can not be handled solely by H\"{o}lder's inequality, primarily because $(2_{*})'$ is larger than $2.$ 

To tackle $S_{1},$ we reformulate it in the following manner: $$S_{1}=\mbox{Re}\;\int_{\rdm\cap (-\infty, t_{0})}\eta_{I}\eta_{Q}^{2}f\cdot \overline{\hat{v}\eta_{I}}=\mbox{Re}\;\int_{\rdm\cap (-\infty, t_{0})}\eta_{I}\eta_{Q}^{2}f\cdot \overline{\omega},$$ where $\omega:=\hat{v}\eta_{I}.$ Set $a_{Q}=\l(\int_{\rz}\eta_{Q}^{2}\r)^{-1}.$ From $\partial_{t}v-\mbox{div}((S(t, x)+D(t, x)-\fint_{2Q}D)\nabla v)=f+\mbox{div} F\; \mbox{in}\;\Omega$ and the fact that $\supp \eta_{Q} \subset 2Q,$ one can easily check that \footnote{ To see this point, one can consult the estimate (2.9) in \cite[Lemma 2.3]{SS}.}, for a.e. $t\in 4I,$ \begin{equation}\label{eq: 4.3}\partial_{t}\tilde{v}(t)=a_{Q}\int_{\rdm} -\l(S(t, x)+D(t, x)-\fint_{2Q}D\r)\nabla v\cdot \nabla \eta_{Q}^{2}+f\eta_{Q}^{2}-F\cdot \nabla \eta_{Q}^{2}.\end{equation} When demonstrating \eqref{eq: 4.3}, the antisymmetry of $D\in L^{\infty}(BMO)$ plays a pivotal role once more because it ensures that $\int_{\rz}D\nabla v\cdot\nabla \eta_{Q}^{2}\in L^{1}(4I)$ by Proposition \ref{pro: 2.1}. Furthermore, a straightforward computation gives that $\omega$ is a weak solution to $$\partial_{t}u-\mbox{div}((S(t, x)+D(t, x)-\fint_{2Q}D)\nabla u)=\tilde{f}+\mbox{div} (F\eta_{I})\; \mbox{in}\;\Omega$$ with \begin{equation*}\begin{aligned}
\tilde{f}:=\eta_{I}\l(f-a_{Q}\int_{\rz}f\eta_{Q}^{2}(x)dx\r)&+a_{Q}\eta_{I}\int_{\rz}\l(S(t, x)+D(t, x)-\fint_{2Q}D\r)\nabla v\cdot \nabla \eta_{Q}^{2} dx\\
&-a_{Q}\eta_{I}\int_{\rz}F\nabla \eta_{Q}^{2}+\int_{\rz}(v(t, x)-\tilde{v}(t))\partial_{t}\eta_{I}\\[4pt]
&:=K_{1}+K_{2},\end{aligned}\end{equation*}  where \begin{equation*}\begin{aligned}K_{1}:=\eta_{I}\l(f-a_{Q}\int_{\rz}f\eta_{Q}^{2}(x)dx\r)&+a_{Q}\eta_{I}\int_{\rz}\l(S(t, x)+D(t, x)-\fint_{2Q}D\r)\nabla v\cdot \nabla \eta_{Q}^{2} dx\\
&-a_{Q}\eta_{I}\int_{\rz}F\nabla \eta_{Q}^{2}.\end{aligned}\end{equation*} Recall that $1<2_{*}<2$ and $f\in L_{loc}^{2_{*}}(\Omega, \cc^{m}).$ By the choice of $\eta,$ using crudely H\"{o}lder inequality, one can prove $$\l(\int_{16I\times 4Q}|K_{1}|^{2_{*}}\r)^{1/2_{*}}\lesssim \l(\int_{16I\times 4Q}|f|^{2_{*}}\r)^{1/2_{*}}+\l(\int_{16I\times 4Q}|\nabla v|^{2}\r)^{1/2}+\l(\int_{16I\times 4Q}|F|^{2}\r)^{1/2}.$$ On the other hand, it follows from the definition of $\tilde{v}$ and a variant of Poincar$\acute{e}$-Sobolev inequality \cite[Lemma 8.3.1]{AD} that \begin{equation*}\begin{aligned}
\l(\int_{16I\times 4Q}|K_{2}|^{2_{*}}\r)^{1/2_{*}}&\lesssim \l(\int_{16I\times 4Q}|v-a_{Q}\int_{\rz}v\eta_{Q}^{2}|^{2}\r)^{1/2}\\
&\lesssim \l(\int_{16I\times 4Q}|v-a_{2Q}\int_{\rz}v\eta_{2Q}^{2}|^{2}\r)^{1/2}\\
&\quad\quad+\l(\int_{16I}\bigg|a_{Q}\int_{\rz}\l(v-\frac{\int_{\rz}v\eta_{2Q}^{2}}{\int_{\rz}\eta_{2Q}^{2}}\r)\eta_{Q}^{2}\bigg|^{2}\r)^{1/2}\\
&\lesssim \l(\int_{16I\times 4Q}|\nabla v|^{2}\r)^{1/2},\quad\quad\end{aligned}\end{equation*} moreover, by the same argument, $$\l(\int_{16I\times 4Q}|\omega|^{2}\r)^{1/2}\lesssim\l(\int_{16I\times 4Q}|\nabla v|^{2}\r)^{1/2}.$$

Equipped with the above estimates, and applying Proposition \ref{pro: 3.4}, we get \begin{align*}
|S_{1}|&\lesssim \l(\int_{16I\times 4Q}|f|^{2_{*}}\r)^{1/2_{*}}\l(\int_{16I\times 4Q}|\omega|^{2^{*}}\r)^{1/2^{*}}\\[4pt]
&\lesssim  \l(\int_{16I\times 4Q}|f|^{2_{*}}\r)^{1/2_{*}}\l(\l(\int_{16I\times 4Q}|\omega|^{2}\r)^{1/2}+\l(\int_{16I\times 4Q}|F\eta_{I}|^{2}\r)^{1/2}+\l(\int_{16I\times 4Q}|\tilde{f}|^{2_{*}}\r)^{1/2_{*}}\r)\\[4pt]
&\lesssim \l(\int_{16I\times 4Q}|f|^{2_{*}}\r)^{1/2_{*}}\l(\l(\int_{16I\times 4Q}|\nabla v|^{2}\r)^{1/2}+\l(\int_{16I\times 4Q}|F|^{2}\r)^{1/2}+\l(\int_{16I\times 4Q}|f|^{2_{*}}\r)^{1/2_{*}}\r)\\[4pt]
&\lesssim \l(\int_{16I\times 4Q}|f|^{2_{*}}\r)^{2/2_{*}}+\delta\l(\int_{16I\times 4Q}|\nabla v|^{2}\r)^{2/2}+\int_{16I\times 4Q}|F|^{2}.
\end{align*} Clearly, by \eqref{eq: 1.2}, $$\mbox{Re}\;\int_{\rdm\cap (-\infty, t_{0})}\eta_{I}^{2}S\nabla(\eta_{Q}\hat{v})\overline{\nabla(\eta_{Q}\hat{v})}\geq \lambda\int_{\rdm\cap (-\infty, t_{0})}|\nabla (\hat{v}\eta)|^{2}-\kappa |(\hat{v}\eta)|^{2}.$$ Then, letting $\delta$ sufficiently small, we arrive at \begin{equation}\label{eq: 4.4}
\begin{aligned}
\sup_{-2<t_{0}<0}\int_{Q}|\hat{v}(t_{0}, x)|^{2}dx+\int_{I\times Q}|\nabla v|^{2}&\lesssim \ez\l(\int_{16I\times 4Q}|\nabla v|^{2}\r)^{2/2}+\l(\int_{16I\times 4Q}|f|^{2_{*}}\r)^{2/2_{*}}\\
&\quad\quad\quad+\int_{16I\times 4Q}|F|^{2}+ \int_{-4}^{0}\l(\int_{2Q}|\hat{v}|^{\frac{2s}{2-s}}\r)^{\frac{2-s}{s}}.
\end{aligned} \end{equation}
Confining $s\leq \frac{n}{n-1},$ and using \cite[Lemma 8.3.1]{AD} again, the last term on the right hand side of \eqref{eq: 4.4} can be bounded by $$\l(\int_{3Q}|\hat{v}|^{\frac{2s}{2-s}}\r)^{\frac{2-s}{s}}\lesssim \int_{4Q}|\nabla v|^{2},$$ by which, also \eqref{eq: 4.4}, we obtain \begin{equation}\label{eq: 4.5}\sup_{-2<t_{0}<0}\int_{Q}|\hat{v}(t_{0}, x)|^{2}dx\lesssim \l(\int_{16I\times 4Q}|\nabla v|^{2}\r)^{2/2}+\l(\int_{16I\times 4Q}|f|^{2_{*}}\r)^{2/2_{*}}+\int_{16I\times 4Q}|F|^{2}.\end{equation}


By further confining $1<s<\frac{2n}{2n-1}<\frac{n}{n-1},$ we observe that for $n=2, 3, ...,$ $\frac{2n}{2n-1}<\frac{4}{3}<\frac{4n}{3n-1}<2.$ Additionally, it is not difficult to determine $0<\lambda, \mu<1$ and $1<r<2$ satisfying $$\frac{2s}{2-s}=2\lambda+\frac{nr}{n-r}\mu,\quad \lambda+\mu=1, \quad\frac{nr}{n-r}\mu\frac{2-s}{s}=1.$$ Going back to \eqref{eq: 4.4}, we can derive \begin{equation}\label{eq: 4.6}\begin{aligned}
\int_{I\times Q}|\nabla v|^{2}&\lesssim \|f\|_{L^{2_{*}}( 16I\times 4Q)}^{2}+\l(\int_{16I\times 4Q}|F|^{2}\r)^{2/2}+\int_{-4}^{0}\l(\int_{2Q}|\hat{v}|^{\frac{2s}{2-s}}\r)^{\frac{2-s}{s}}\\
&\quad\quad+\ez\int_{16I\times 4Q}|\nabla v|^{2}\\
&\approx \|f\|_{L^{2_{*}}( 16I\times 4Q)}^{2}+\l(\int_{16I\times 4Q}|F|^{2}\r)^{2/2}+\int_{-4}^{0}\l(\int_{2Q}|\hat{v}|^{2\lambda+\frac{nr}{n-r}\mu}\r)^{\frac{2-s}{s}}\\
&\quad\quad +\ez\int_{16I\times 4Q}|\nabla v|^{2}\\
&\lesssim \|f\|_{L^{2_{*}}( 16I\times 4Q)}^{2}+\l(\int_{16I\times 4Q}|F|^{2}\r)^{2/2}+\int_{-4}^{0}\l(\int_{2Q}|\hat{v}|^{2}\r)^{\frac{1}{2}}\l(\int_{2Q}|\hat{v}|^{\frac{nr}{n-r}}\r)^{\frac{n-r}{nr}}\\
&\quad\quad+\ez\int_{16I\times 4Q}|\nabla v|^{2}\\
&\lesssim \|f\|_{L^{2_{*}}( 16I\times 4Q)}^{2}+\l(\int_{16I\times 4Q}|F|^{2}\r)^{2/2}++\ez\int_{16I\times 4Q}|\nabla v|^{2}\end{aligned} \end{equation}
\begin{equation*}\begin{aligned}
&\quad\quad+\sup_{-4<t_{0}<0}\l(\int_{2Q}|\hat{v}|^{2}\r)^{1/2}\l(\int_{4I\times 2Q}|\nabla v|^{r}\r)^{1/r}.
\end{aligned} \end{equation*}\\
Replacing $I\times Q$ with $4I\times 2Q$ in \eqref{eq: 4.5}, it follows that \begin{equation}\label{eq: 4.7}\begin{aligned}
\sup_{-4<t_{0}<0}\int_{2Q}|\hat{v}(t, x)|^{2}&=\sup_{-4<t_{0}<0}\int_{2Q}\bigg|v(t, x)-\frac{\int_{\rz}v(t, x)\eta_{Q}(x)^{2}dx}{\int_{\rz}\eta_{Q}(x)^{2}}\bigg|^{2}\\
&\lesssim\sup_{-4<t_{0}<0}\int_{2Q}\bigg|v(t, x)-\frac{\int_{\rz}v(t, x)\eta_{2Q}(x)^{2}dx}{\int_{\rz}\eta_{2Q}(x)^{2}}\bigg|^{2}\\
&\quad+\sup_{-4<t_{0}<0}\bigg|\frac{\int_{\rz}\l(v(t, x)-\frac{\int_{\rz}v(t, x)\eta_{2Q}(x)^{2}dx}{\int_{\rz}\eta_{2Q}(x)^{2}}\r)\eta_{Q}(x)^{2}dx}{\int_{\rz}\eta_{Q}(x)^{2}}\bigg|^{2}\\
&\lesssim \sup_{-4<t_{0}<0}2\int_{2Q}\bigg|v(t, x)-\frac{\int_{\rz}v(t, x)\eta_{2Q}(x)^{2}dx}{\int_{\rz}\eta_{2Q}(x)^{2}}\bigg|^{2}\\
&\lesssim \|f\|_{L^{2_{*}}( 16I\times 4Q)}^{2}+\l(\int_{16I\times 3Q}|F|^{2}\r)^{2/2}+\int_{16I\times 4Q}|\nabla v|^{2}.
\end{aligned} \end{equation}
Combining \eqref{eq: 4.6} and \eqref{eq: 4.7}, also applying Young's inequality, we instantly conclude \eqref{eq: 4.2}, equivalently, \begin{equation}\label{eq: 4.8}\begin{aligned}\l(\fint_{I\times Q}|\nabla v|^{2}\r)^{1/2}&\lesssim \l(\fint_{16I\times 4Q}|f|^{2_{*}}\r)^{1/2_{*}}+\l(\fint_{16I\times 4Q}|F|^{2}\r)^{1/2}\\
&\quad\quad+\l(\fint_{16I\times 4Q}|\nabla v|^{r}\r)^{1/r}+\ez \l(\fint_{16I\times 4Q}|\nabla v|^{2}\r)^{1/2}.\end{aligned} \end{equation} Using \eqref{eq: 4.8} and employing a standard covering argument, followed by scaling, we eventually deduce that \begin{equation}\label{eq: 4.9}\begin{aligned}\l(\fint_{I_{R}\times Q_{R}}|\nabla v|^{2}\r)^{1/2}&\lesssim R\l(\fint_{4I_{R}\times 2Q_{R}}|f|^{2_{*}}\r)^{1/2_{*}}+\l(\fint_{4I_{R}\times 2Q_{R}}|F|^{2}\r)^{1/2}\\
&\quad\quad+\l(\fint_{4I_{R}\times 2Q_{R}}|\nabla v|^{r}\r)^{1/r}+\ez \l(\fint_{4I_{R}\times 2Q_{R}}|\nabla v|^{2}\r)^{1/2},\end{aligned} \end{equation} where $I_{R}:=(t_{0}-R^{2}, t_{0})$ and $Q_{R}:=B_{R}(x_{0})$ such that $4I_{R}\times 2Q_{R}\subset \Omega$ for any $(t_{0}, x_{0})\in \Omega.$ 

In order to obtain \eqref{eq: 4.1}, we need to demonstrate that the estimate \eqref{eq: 4.9} makes it possible to apply the improved Gehring lemma, that is Lemma \ref{lemma: 5.1} in the appendix. For this goal, we introduce $$g_{r}:=|\nabla v|^{r},\quad f_{r}:=|f|^{r},\quad F_{r}:=|F|^{r}.$$ Then \eqref{eq: 4.9} becomes \begin{equation}\label{eq: 4.10}\begin{aligned}\l(\fint_{I\times Q}g_{r}^{\frac{2}{r}}\r)^{r/2}&\lesssim R^{r}\l(\fint_{16I\times 4Q}f_{r}^{\frac{2_{*}}{r}}\r)^{r/2_{*}}+\l(\fint_{16I\times 4Q}F_{r}^{\frac{2}{r}}\r)^{r/2}+\fint_{16I\times 4Q}g_{r}\\
&\quad\quad+\ez^{r} \l(\fint_{16I\times 4Q}g_{r}^{\frac{2}{r}}\r)^{r/2}.\end{aligned} \end{equation} Set $q:=\frac{2}{r},$ $s=\frac{2_{*}}{r}$ and $\beta:=r=(\frac{1}{s}-\frac{1}{q})(n+2).$ Bear in mind that we allow for $0<s<1$ in Lemma \ref{lemma: 5.1}. Note that $q>1,$ and from \eqref{eq: 4.10}, the conditions of Lemma \ref{lemma: 5.1} are satisfied. As a result, there exists a $p_{r}>q$ small enough such that \begin{equation}\label{eq: 4.11}\begin{aligned}\l(\fint_{I\times Q}g_{r}^{p_{r}}\r)^{1/p_{r}}&\lesssim R^{r}\l(\fint_{16I\times 4Q}f_{r}^{(p_{r})_{*, \beta}}\r)^{1/(p_{r})_{*, \beta}}+\l(\fint_{16I\times 4Q}F_{r}^{p_{r}}\r)^{1/p_{r}}\\
&\quad\quad\quad+\l(\fint_{16I\times 4Q}g_{r}^{q}\r)^{1/q}.\end{aligned} \end{equation} Recall that $(p_{r})_{*, \beta}=\frac{(n+2)p_{r}}{n+2+\beta p_{r}}.$ Letting $p:=rp_{r},$ it follows that $p>2$ is sufficiently small and $(p_{r})_{*, \beta}=\frac{p_{*}}{r}.$ Then, \eqref{eq: 4.11} becomes \begin{equation}\label{eq: 4.12}\begin{aligned}\l(\fint_{I\times Q}|\nabla v|^{p}\r)^{1/p}&\lesssim R \l(\fint_{16I\times 4Q}|f|^{p_{*}}\r)^{1/p_{*}}+\l(\fint_{16I\times 4Q}|F|^{p}\r)^{1/p}\\&\quad\quad+\l(\fint_{16I\times 4Q}|\nabla v|^{2}\r)^{1/2}.\end{aligned} \end{equation} Finally,  \eqref{eq: 4.12} and Theorem \ref{thm: 3.5} lead to \eqref{eq: 4.1}, thus concluding the proof.  

\end{proof}

\section{Appendix}

As we focus on local estimates, defining the coefficient matrices $S$ and $D$ globally are not required. Indeed, the following extension lemma holds. 

\begin{lemma}\label{lemma: 5.2}\; For some open interval $I_{0}\subset \rr$ and open ball (hence NTA domain) $Q_{0}\subset \rdm.$ 
\begin{equation*}
\begin{aligned} 
&(i)\quad \mbox{Assume}\; S\in L^{\infty}(\rr^{n+1}, \cc^{mn})\;\mbox{satisfies the G$\mathring{a}$rding inequality (uniformly in}\;t\in I_{0}):\\
&\quad\quad \mbox{Re} \;\int_{Q_{0}}\;S(t,x)\nabla u(t, x)\cdot \overline{\nabla u(t, x)}\geq \lambda\int_{Q_{0}}|\nabla u|^{2}-\kappa \int_{Q_{0}}| u|^{2},\quad \forall\;u\in W^{1,2}_{0}(Q_{0}),\\
&\mbox{then, given $\ez\in (0, 1),$ there exists a $\tilde{S}$ with $\tilde{S}=S$ on $(1-\ez)^{2}I_{0}\times (1-\ez)Q_{0}$ that verifies \eqref{eq: 1.2},}\\
&\mbox{where the ellipticity constant for $\tilde{S}$ are possibly different and may depend on $\ez.$ }\\
&(ii)\quad \mbox{Assume $D\in L^{\infty}(I_{0}; BMO(Q_{0}))$ is real and anti-symmetric, then there exists another  real}\\
&\mbox{ and anti-symmetric matrix $\tilde{D},$ belonging to $L^{\infty}(\rr; BMO(\rz)),$ such that }\\
&\quad\quad\quad\quad\quad\quad\quad \tilde{D}|_{I_{0}\times Q_{0}}=D\quad\mbox{and}\quad \|\tilde{D}\|_{L^{\infty}(\rr; BMO(\rz))} \lesssim \|D\|_{L^{\infty}(I_{0}; BMO(Q_{0}))}.\\
\end{aligned} 
\end{equation*}
\end{lemma}
 
{\it Proof.}\quad $(i)$ is proved in \cite[Lemma A.1]{ABES}. For any fixed $t\in I_{0},$ since $D(t, \cdot)\in BMO(Q_{0}),$ it follows from \cite{J} that there exists $\overline{D}(t, \cdot)\in BMO(\rz)$ such that $$\overline{D}(t, \cdot)|_{ Q_{0}}=D(t, \cdot)\quad\mbox{and}\quad \|\overline{D}(t, \cdot)\|_{BMO(\rz)} \leq c\|D(t, \cdot)\|_{ BMO(Q_{0})},$$ where the constant $c$ depends only on the domain and dimension. Set $$\tilde{D}:=\l\{
\begin{aligned}
&\;\overline{D}(t, \cdot),\quad \mbox{if}\; t\in I_{0}\\
&\;0,\quad\quad\quad \; \mbox{if} \;t\in \rr \setminus I_{0}.\end{aligned}\r.$$ Clearly, $\tilde{D}$ meets the requirements of $(ii).$ 

\hfill$\Box$

The classical Gehring lemma (\cite[Chapter \;V \;Proposition\; 1.1]{GM}) is useful in proving H\"{o}lder type estimates for elliptic and parabolic equations (systems). In this section, we develop an upgraded version (Lemma \ref{lemma: 5.1}), emphasizing that the exponent $\frac{(n+2)p}{n+2+\beta p}$ is strictly less than $\frac{ps}{q},$ as provided by \cite[Chapter \;V \;Proposition\; 1.1]{GM}. For a detailed derivation of the exponent $\frac{ps}{q},$ we refer to \cite[Chapter \;V\; Theorem 2.2]{GM} (or \cite[Theorem 2.1]{ABESa}). As a consequence, the exponent governing the sequence of functions $\{f_{i}\}$ in \cite[Chapter \;V \;Theorem\; 2.2]{GM} can be improved. Such a result also has impact on partial regularity of nonlinear systems on \cite{C, GG, GS}.

The critical exponent $\frac{(n+2)p}{n+2+\beta p}$ emerges from the non-local Gehring lemma in \cite[Theorem 4.2]{ABESa}, which is formulated for spaces of homogeneous type. However, the proof, more specifically \cite[Theorem 2.1]{ABESa}, appears to be inapplicable to the local case with an additionally small tail, like the $\theta-$dependent term in \eqref{eq: 5.1}. To overcome the difficulty, we revisit the method outlined in \cite[Chapter \;V \;Proposition\; 1.1]{GM}, taking full advantage of the weak and strong boundedness properties of the fractional maximal function on $L^{p}(1\leq p<\infty)$. As shown in Theorem \ref{thm: 4.2}, Lemma \ref{lemma: 5.1} is more effective than the classical approach for enhancing the regularity of solutions to linear (or nonlinear) elliptic and parabolic equations, including their systems, with inhomogeneous low order terms. 

\begin{lemma}\label{lemma: 5.1}\;  Given a parabolic cube $Q\subset \rdm.$ 
Let $g, h\in L^{q}(Q)(q>1)$ and $f^{s}\in L^{1}(Q) (0<s<q)$ be three nonnegative functions, and set $\beta:=(n+2)(\frac{1}{s}-\frac{1}{q}).$ Suppose that \begin{equation}\label{eq: 5.1}\fint_{Q_{r}(x_{0})}g^{q}\leq a\l(\fint_{Q_{2r}(x_{0})}g\r)^{q}+\fint_{Q_{2r}(x_{0})}h^{q}+r^{\beta q}\l(\fint_{Q_{2r}(x_{0})}f^{s}\r)^{\frac{q}{s}}+\theta \fint_{Q_{2r}(x_{0})}g^{q}
\end{equation} holds for each $x_{0}\in Q$ with each $r<\frac{1}{2}\min(dist(x_{0}, \partial Q), r_{0}),$ where $r_{0}, a, \theta$ are constants with $a>1,$ $r_{0}>0,$ $0\leq \theta <1.$ Then $g\in L^{p}_{loc}(Q)$ for $p\in [q, q+\ez)$ and \begin{equation}\label{eq: 5.2}\l(\fint_{Q_{r}}g^{p}\r)^{1/p}\leq b\l(\l(\fint_{Q_{2r}}g^{q}\r)^{1/q}+\l(\fint_{Q_{2r}}h^{p}\r)^{1/p}+r^{\beta}\l(\fint_{Q_{2r}}f^{p_{*, \beta}}\r)^{1/p_{*, \beta}}\r)\end{equation} for $Q_{2r}\subset Q, r<r_{0},$ where $p_{*, \beta}:=\frac{(n+2)p}{n+2+\beta p}$ and the two constants $b$ and $\ez$ depend only on $a, \theta, q, n$ and $s$.
\end{lemma}
 
{\it Proof.}\quad The case $\theta =0$ is indeed covered in \cite[Theorem 4.2]{ABESa}, so we only need to address the case $0<\theta <1.$ Note that the following method also applies to $\theta =0$. 

In line with the procedure used in the proof of \cite[Chapter \;V \;Proposition\; 1.1]{GM}, our initial task is to build that if almost every on $Q:=Q_{3/2}(0)$ \begin{equation}\label{eq: 5.3}M_{\frac{1}{2}d(x)}(g^{q})(x)\leq a M^{q}(g)(x)+M(h^{q})+M^{\beta s}(f^{s})^{\frac{q}{s}}(x)+\theta M(g^{q})(x),\end{equation} where $d(x):=dist(x, \partial Q),$ $M$ denotes the centered Hardy-Littlewood maximal function and $$M_{r_{0}}(u)(x):=\sup_{r(B)<r_{0}}\fint_{B_{r}(x)}|u|, \quad M^{\alpha}(u)(x):=\sup_{B\ni x}r(B)^{\alpha}\fint_{B}|u|,$$  then \begin{equation}\label{eq: 5.4}\l(\fint_{Q_{1/2}(0)}g^{p}\r)^{1/p}\lesssim \l(\fint_{Q_{3/2}(0)}g^{q}\r)^{1/q}+\l(\fint_{Q_{3/2}(0)}h^{p}\r)^{1/p}+\l(\fint_{Q_{3/2}(0)}f^{p_{*, \beta}}\r)^{1/p_{*, \beta}}.\end{equation} Here, the implicit constant is determined only by $a, \theta, q, n, \tau$ and $s$. We explicitly assume $g=h=f=0$ in $\rdm\setminus Q$ in the subsequent analysis.

We now split $Q$ into $$Q=\cup_{k}C_{k},\quad C_{0}:=\{x=(t, x_{1}, ...x_{n})\in R^{n+1}: |x_{j}|<1/2, \;|t|<1/4\}$$ and $$C_{k}:=\{x\in Q: 2^{-k}<d(x)\leq 2^{-k+1}\}, \quad k\geq 1.$$ Set $\sigma>100D^{1/2}$ with $D:=n+2.$ We also define $\alpha_{k}:=\l(\sigma^{D}2^{kD}\r)^{1/q},$ $E(u, \lambda):=\{x\in Q: u(x)>\lambda\}$, $$G(x):=\frac{\theta^{1/2} g(x)}{\|h\|_{L^{q}(Q)}+\|g\|_{L^{q}(Q)}+\|f\|_{L^{s}(Q)}}, \quad \tilde{G}(x)=\alpha_{k}^{-1}G(x)\;\mbox{in}\; C_{k},$$ $$F(x):=\frac{\theta^{1/2} f(x)}{\|h\|_{L^{q}(Q)}+\|g\|_{L^{q}(Q)}+\|f\|_{L^{s}(Q)}}, $$ and $$H(x):=\frac{\theta^{1/2} h(x)}{\|h\|_{L^{q}(Q)}+\|g\|_{L^{q}(Q)}+\|f\|_{L^{s}(Q)}}.$$ With these preparations in place, the inequality \eqref{eq: 5.3} transforms into \begin{equation}\label{eq: 5.5}M_{\frac{1}{2}d(x)}(G^{q})(x)\leq c(n, a)M^{q}[G+M^{\beta s}(F^{s})^{1/s}+M(H^{q})^{1/q}](x)+\theta M(G^{q}(x)) \quad\;\mbox{a.e.}.\end{equation} In the above derivation, the well-known fact that $u(x)\leq M(u)(x)$ was exploited. In what follows, for the sake of brevity, we let $$U:=M^{\beta s}(F^{s})^{1/s}+M(H^{q})^{1/q}.$$ 

Let $\lambda$ be a fixed real number in $[1, \infty).$ Define $$\mu:=\delta\cdot \lambda,$$ where $\delta>1$ is a constant, to be determined later. According to the definition of $G,$ $$\fint_{Q}G^{q}\leq \theta^{1/2}<\mu^{q}.$$ Consequently, applying the Calder\'{o}n-Zygmund decomposition-specifically, a parabolic version of \cite[Chapter \;V \;Lemma\; 1.3]{GM}-we obtain a disjoint sequence of sub-cubes $\{Q_{k}^{j}\}$ of $Q$ such that \begin{equation}\label{eq: 5.6}Q_{k}^{j}\subset C_{k},\; \forall j, k\in \nn, \quad diam(Q_{k}^{j})<\frac{1}{2}dist(Q_{k}^{j}, \partial Q),\end{equation}
\begin{equation}\label{eq: 5.7} (\alpha_{k}\mu)^{q}<\fint_{Q_{k}^{j}}G^{q}\leq \sigma^{n}(\alpha_{k}\mu)^{q},\end{equation} and $$G\leq c(n, q)\alpha_{k}\mu\quad \mbox{in}\; C_{k}\setminus \cup_{j}Q_{k}^{j}. $$
The last two inequalities are plainly equivalent to $$\mu^{q}<\fint_{Q_{k}^{j}}\tilde{G}^{q}\leq \sigma^{n}\mu^{q} \quad\mbox{and}\quad \tilde{G} \leq c(n, q)\mu\quad \mbox{in}\; C_{k}\setminus \cup_{j}Q_{k}^{j}.$$ Since $|E(\tilde{G}, c(n, q)\mu)\setminus \cup_{j, k}Q_{k}^{j}|=0,$ then \begin{equation}\label{eq: 5.8}\int_{E(\tilde{G}, c(n, q)\mu)}\tilde{G}^{q}\leq \sum_{j, k}\int_{Q_{k}^{j}}\tilde{G}^{q}\leq \sigma^{n}\mu^{q}\sum_{j, k}|Q_{k}^{j}|.\end{equation} If we set $\tilde{Q}:=\{x\in Q: \eqref{eq: 5.5} \; \mbox{is true}\}$ and $\hat{Q}:=\cup_{j, k}Q_{k}^{j},$ then $|Q|=|\tilde{Q}|$ and $|Q_{k}^{j}|=|Q_{k}^{j}\cap \tilde{Q}|.$ In the sequel we seek to control $|\hat{Q}|=\sum_{j, k}|Q_{k}^{j}|=\sum_{j, k}|Q_{k}^{j}\cap \tilde{Q}|. $ 

To the end, we first focus on annihilating the term $\theta M(G^{q}(x)).$ Given any $x\in Q_{k}^{j},$ by the definition of $M_{\frac{1}{2}d(x)}$ and \eqref{eq: 5.6}-\eqref{eq: 5.7}, we see \begin{equation}\label{eq: 5.9}(\alpha_{k}\mu)^{q}<\fint_{Q_{k}^{j}}G^{q}\leq c(n)M_{\frac{1}{2}d(x)}(G^{q})(x)\leq c(n)M(G^{q})(x).\end{equation}This implies that there exists a parabolic ball $B$ centered at $x$ such that \begin{equation}\label{eq: 5.10}c(n)^{-1}(\alpha_{k}\mu)^{q}\theta^{1/2}\leq \theta^{1/2}M(G^{q})(x) <\fint_{B}G^{q}.\end{equation} Furthermore, keeping in mind that $\mu>1,$ the latter expression readily leads to the conclusion that $$|B|\leq \frac{c(n)}{\theta^{1/2}}\frac{1}{\alpha_{k}^{q}}\int_{B}G^{q}<\frac{c(n)}{\alpha_{k}^{q}},$$ from which we get $$r(B)<\frac{d(x)}{2}.$$ This, together with\eqref{eq: 5.10}, in turn implies that  $$M(G^{q})(x)\leq \frac{1}{\theta^{1/2}}M_{\frac{1}{2}d(x)}(G^{q})(x).$$ Hence, if $x\in Q_{k}^{j}\cap \tilde{Q},$ \eqref{eq: 5.5} indicates that \begin{equation}\label{eq: 5.11}M_{\frac{1}{2}d(x)}(G^{q})(x)\leq \frac{c(n, a)}{1-\theta^{1/2}}M^{q}[G+U](x).\end{equation} Plugging  \eqref{eq: 5.11} into  \eqref{eq: 5.9} we acheive\begin{equation*}\begin{aligned}(\alpha_{k}\mu)^{q}<\fint_{Q_{k}^{j}}G^{q}&\leq c(n)M_{\frac{1}{2}d(x)}(G^{q})(x)\\
&\leq c(n)M(G^{q})(x)\leq \frac{c(n, a)}{1-\theta^{1/2}}M^{q}[G+U](x),\end{aligned}\end{equation*}
Thus, we can find a parabolic ball $B_{x},$ centered at $x,$ satisfying the following property: \begin{equation}\label{eq: 5.12}(\alpha_{k}\mu)^{q}\leq \frac{c(n, a)}{1-\theta^{1/2}}\l(\fint_{B_{x}} (G+U)\r)^{q}. \end{equation} Moreover, invoking the definition of $G,$ we deduce that $$\l(\fint_{B_{x}}G\r)^{q}\leq \frac{c(n, q)}{|B_{x}|}\int_{Q}G^{q}\leq \frac{c(n, q, \theta)}{|B_{x}|}.$$ Regarding the average of $U$, one can prove by the $L^{1, \infty}(\rdm)\to L^{1}(\rdm)$ boundedness of $M$ and the H\"{o}lder inequality in Lorentz spaces that  \begin{align*}\l(\fint_{B_{x}} (M(H^{q})^{1/q}\r)^{q}&\leq \frac{c(n, q)}{|B_{x}|} \|M(H^{q})\|_{L^{1, \infty}(B_{x})}\\
&\leq \frac{c(n, q)}{|B_{x}|}\int_{Q}H^{q}\leq \frac{c(n, q, \theta)}{|B_{x}|}. \end{align*} Recall that $\beta=D(\frac{1}{s}-\frac{1}{q})$ and $q>s.$ Then $\frac{q}{s}=\frac{D}{D-\beta s},$ and the weak $(1, \frac{q}{s})$ boundedness of $M^{\beta s}$ is true. Utilizing this along with the H\"{o}lder inequality in Lorentz spaces, we arrive at \begin{align*}
\l(\fint_{B_{x}} (M^{\beta s}(F^{s})^{1/s}\r)^{q}&\leq \frac{c(n, q)}{|B_{x}|} \|M^{\beta s}(F^{s})^{1/s}\|_{L^{q, \infty}(B_{x})}^{q}\\[4pt]
&\leq \frac{c(n, q)}{|B_{x}|} \|M^{\beta s}(F^{s})\|_{L^{q/s, \infty}(\rdm)}^{q/s}\\[4pt]
&\leq \frac{c(n, q, s)}{|B_{x}|}\l(\int_{Q}|F|^{s}\r)^{q/s}\leq \frac{c(n, q, s, \theta)}{|B_{x}|}.
\end{align*} Summarizing all the estimates above we conclude from \eqref{eq: 5.12} that $$(\alpha_{k}\mu)^{q}\leq \frac{c(n, q, s,\theta, a)}{|B_{x}|},$$ which, by $\lambda>1,$ gives that  $$r(B_{x})\leq \frac{c(n, q, s, \theta, a)}{\delta}2^{-k}.$$ 

If $\delta=\delta(n, q, s, \theta, a)$ is choosen large enough, the latter estiamte guarantees that the ball $B$ at most intersects with three annuli $C_{k-1}, C_{k}, C_{k+1}.$ Bear in mind  that $\alpha_{k+1}=2^{D/q}\alpha_{k}.$ Returning to \eqref{eq: 5.12} we find that $$\mu\leq \l(\frac{c(n, a)}{1-\theta^{1/2}}\r)^{1/q} \frac{1}{\alpha_{k}}\fint_{B_{x}} (G+U)\leq c(n, a, \theta, q)\fint_{B_{x}}(\tilde{G}+\tilde{U}),$$ where \begin{equation}\label{eq: 5.13}\tilde{U}(x):=\alpha_{k}^{-1}\l(M^{\beta s}(F^{s})^{1/s}(x)+M(H^{q})^{1/q}(x)\r)\quad \mbox{in}\; C_{k}.\end{equation} From \eqref{eq: 5.11} and the definition of $\mu,$ it is apparant that $$\delta \lambda |B_{x}| \leq c(n, a, \theta, q)\l[\int_{B_{x}\cap E(\tilde{G}, \lambda)}\tilde{G}+\int_{B_{x}\cap E(\tilde{U}, \lambda)}\tilde{U}+2\lambda |B_{x}|\r].$$ Then, letting $\delta>2c(n, a, \theta, q),$ we see \begin{equation}\label{eq: 5.14}|B_{x}|\leq \frac{c(n, a, \theta, q)}{\lambda}\l[\int_{B_{x}\cap E(\tilde{G}, \lambda)}\tilde{G}+\int_{B_{x}\cap E(\tilde{U}, \lambda)}\tilde{U}\r].\end{equation} As $\hat{Q}\cap \tilde{Q}\subset \cup_{x\in \hat{Q}\cap \tilde{Q}}B_{x},$ the Vitali's covering lemma (a parabolic version of \cite[Chapter \;V \;Lemma\; 1.1]{GM}) yields that there exists a disjoint subsequence $\{B_{l}\}$ such that $$|\hat{Q}\cap \tilde{Q}|\leq c(n) \sum_{l}|B_{l}|.$$ Recalling \eqref{eq: 5.8}, also using \eqref{eq: 5.14}, we get \begin{align*}
\int_{E(\tilde{G}, c(n, q)\mu)}\tilde{G}^{q}&\leq c(n, q)\lambda^{q}\sum_{l}|B_{l}|\\[4pt]
&\leq c(n, q, a, \theta)\lambda^{q-1}\sum_{l}\l[\int_{B_{l}\cap E(\tilde{G}, \lambda)}\tilde{G}+\int_{B_{l}\cap E(\tilde{U}, \lambda)}\tilde{U}\r]\\[4pt]
&\leq c(n, q, a, \theta)\lambda^{q-1}\l[\int_{ E(\tilde{G}, \lambda)}\tilde{G}+\int_{ E(\tilde{U}, \lambda)}\tilde{U}\r].\end{align*} Combining this with $$\int_{E(\tilde{G}, \lambda)\setminus E(\tilde{G}, c(n, q)\delta \lambda)}\tilde{G}^{q}\leq c(n, q, a, \theta, s)\lambda^{q-1}\int_{ E(\tilde{G}, \lambda)}\tilde{G},$$ we arrive at \begin{equation}\label{eq: 5.15}\int_{E(\tilde{G}, \lambda )}\tilde{G}^{q}\leq c(n, q, a, s, \theta)\lambda^{q-1}\l[\int_{ E(\tilde{G}, \lambda)}\tilde{G}+\int_{E(\tilde{U}, \lambda)}\tilde{U}\r], \;\forall \;\lambda\geq 1.\end{equation}

Next, we derive \eqref{eq: 5.4} from \eqref{eq: 5.15}. For this goal, we set $$v(\lambda):=\int_{ E(\tilde{G}, \lambda)}\tilde{G},\quad V(\lambda):=\int_{E(\tilde{U}, \lambda)}\tilde{U}, \quad \mbox{for}\; \lambda\geq 1.$$ By expressing the $L^{p}$ norm in terms of the distribution function, we have that for any nonnegative function $u$ and $\tau\geq 1,$ \begin{equation}\label{eq: 5.16}\int_{E(u, \lambda)}u^{\tau}=-\int_{\lambda}^{\infty}u^{\tau-1}d\l(\int_{E(u, t)}u\r).\end{equation} Observe that, due to \eqref{eq: 5.15}, the requirements of \cite[Chapter \;V \;Lemma\; 1.2]{GM} are fulfilled (with $h, H, p, q$ replaced by $v, V, p-1, q-1$). Therefore, for $p\in [q, q+\ez),$ $$-\int_{1}^{\infty}\lambda^{p-1}dv(\lambda)\leq c(n, q, a, s, \theta, p)\l[-\int_{1}^{\infty}\lambda^{q-1}dv(\lambda)+-\int_{1}^{\infty}\lambda^{p-1}dV(\lambda)\r],$$ which, in conjunction with \eqref{eq: 5.16}, gives $$\int_{E(\tilde{G}, 1)}\tilde{G}^{p}\leq c(n, q, a, s, \theta, p)\l(\int_{E(\tilde{G}, 1)}\tilde{G}^{q}+\int_{E(\tilde{U}, 1)}\tilde{U}^{p}\r).$$ Moreover, as $p>q,$ $$\int_{\{x\in Q: \tilde{G}\leq 1\}}\tilde{G}^{p}=\int_{\{x\in Q: \tilde{G}\leq 1\}}\tilde{G}^{q}\tilde{G}^{p-q}\leq \int_{\{x\in Q: \tilde{G}\leq 1\}}\tilde{G}^{q}.$$ Gathering the above inequalities, we conclude that  
\begin{equation}\label{eq: 5.17}\int_{Q}\tilde{G}^{p}\leq c(n, q, a, s, \theta, p)\l(\int_{Q}\tilde{G}^{q}+\int_{Q}\tilde{U}^{p}\r).\end{equation}

Having obtained \eqref{eq: 5.17}, we are now in a position to infer \eqref{eq: 5.4}.
Since $\beta=D(\frac{1}{s}-\frac{1}{q}),$ it follows that $\beta s<D(p-s).$ This implies $\gamma:=\frac{p_{*, \beta}}{s}>1$ and $\frac{p}{s}=\frac{\gamma D}{D-\beta s\gamma}.$ As shown in \cite{HKNT}, $M^{\beta s}$ is bounded form $L^{\gamma}(\rdm)$ to $L^{\frac{p}{s}}(\rdm).$ Thus \begin{equation}\label{eq: 5.18}
\int_{Q}M^{\beta s}(F^{s})^{p/s}\leq c(n, q, p, s) \l(\int_{Q}F^{s\gamma}\r)^{\frac{p}{s \gamma}}\approx c(n, q, a, \theta)\l(\int_{Q}F^{p_{*, \beta}}\r)^{\frac{p}{p_{*, \beta}}}.\end{equation} On the other hand, \begin{equation}\label{eq: 5.19}\int_{Q}M(H^{q})^{p/q}\leq c(n, q, p) \int_{Q}H^{p}\end{equation} thanks to the fact that $M$ is bounded from $L^{\gamma}(\rdm)$ to $L^{\gamma}(\rdm)$ for any $1<\gamma<\infty.$ Employing \eqref{eq: 5.13} ($\tilde{U}\leq c(n, q) U(x)$ in $Q$) and \eqref{eq: 5.17}-\eqref{eq: 5.19}, we deduce that \begin{equation}\label{eq: 5.20}\l(\int_{Q}\tilde{G}^{p}\r)\leq c(n, q, a, s, \theta, p)\l[\l(\int_{Q}\tilde{G}^{q}\r)+\l(\int_{Q}F^{p_{*, \beta}}\r)^{p/p_{*, \beta}}+\l(\int_{Q}H^{p}\r)\r].\end{equation} Based on the definitions of $G, F, H$ and $\tilde{G},$ and utilizing Young's inequality, H\"{o}lder's inequality and the inequalities $p_{*, \beta}>s$ and $p>q,$ \eqref{eq: 5.20} leads to the conclusion that \begin{align*}
\int_{Q_{1/2}(0)}|g|^{p}&\lesssim \l(\|h\|_{L^{q}(Q)}+\|g\|_{L^{q}(Q)}+\|f\|_{L^{s}(Q)}\r)^{p-q}\int_{Q}|g|^{q}\\[4pt]
&\quad\quad+\l(\int_{Q}f^{p_{*, \beta}}\r)^{p/p_{*, \beta}}+\int_{Q}h^{p}\\[4pt]
&\lesssim \|h\|_{L^{q}(Q)}^{p-q}\int_{Q}|g|^{q} +\|f\|_{L^{s}(Q)}^{p-q}\int_{Q}|g|^{q} +\l(\int_{Q}|g|^{q}\r)^{p/q}+\l(\int_{Q}f^{p_{*, \beta}}\r)^{p/p_{*, \beta}}+\int_{Q}h^{p}\\[4pt]
&\lesssim \int_{Q}h^{p}+\l(\int_{Q}|g|^{q}\r)^{p/q}+\l(\int_{Q}f^{p_{*, \beta}}\r)^{p/p_{*, \beta}},
\end{align*} which readily yields \eqref{eq: 5.4}.  

At this point, the key remaining step is to establish that the assumption \eqref{eq: 5.1} can be attributed to \eqref{eq: 5.3} by translations and dilatations. In fact, if we assume $Q=Q_{R}(z)$ and set $$g_{R}^{z}(y):=g(Ry+z),\quad f_{R}^{z}(y):=R^{\beta} f(Ry+z), \quad h_{R}^{z}(y):=h(Ry+z), \quad \tau:=\frac{r}{R},\quad \xi:=\frac{x_{0}-z}{R},$$ then, by \eqref{eq: 5.1},  \begin{equation*}\begin{aligned}\l(\fint_{Q_{\tau}(\xi)}(g_{R}^{z})^{q}\r)^{1/q}&\leq c\l(\fint_{Q_{2\tau}(\xi)}g_{R}^{z}\r)+\l(\fint_{Q_{2\tau}(\xi)}(h_{R}^{z})^{q}\r)^{1/q}\\
&\quad\quad+\tau^{\beta }\l(\fint_{Q_{2\tau}(\xi)}(f_{R}^{z})^{s}\r)^{\frac{1}{s}}+\theta^{1/q} \l(\fint_{Q_{2\tau}(\xi)}(g_{R}^{z})^{q}\r)^{1/q}\end{aligned}\end{equation*} for each $\xi\in Q_{1}(0)$ with each $\tau<\frac{1}{2}\min(dist(\xi, \partial Q_{1}(0)), \frac{r_{0}}{R}).$ This implies that $$M_{\frac{1}{2}d(\xi)\wedge \frac{r_{0}}{R}}((g_{R}^{z})^{q})(\xi)\leq a M^{q}(g_{R}^{z})(x)+M((h_{R}^{z})^{q})+M^{\beta s}((f_{R}^{z})^{s})^{\frac{q}{s}}(\xi)+\theta M((g_{R}^{z})^{q})(\xi).$$ We would like to point out that, throughout our applications, the value of $r_{0}$ is always taken to be equal to $l(Q)=R.$ As a result, by \eqref{eq: 5.4}, we eventually achieve $$\l(\fint_{Q_{1/2}(0)}(g_{R}^{z})^{p}\r)^{1/p}\leq b \l(\fint_{Q_{1}(0)}(g_{R}^{z})^{q}\r)^{1/q}+\l(\fint_{Q_{1}(0)}(h_{R}^{z})^{p}\r)^{1/p}+\l(\fint_{Q_{1}(0)}(f_{R}^{z})^{p_{*, \beta}}\r)^{1/p_{*, \beta}}.$$ The inequality is exactly equivalent to \eqref{eq: 5.2}, thus completing the entire proof.

\hfill$\Box$
 
\section*{Availability of data and material}
 Not applicable.
 
 \section*{Competing interests}
 The author declares that they have no competing interests.

\end{document}